\documentclass[reqno]{amsart}  
\theoremstyle{plain}
\newtheorem{theorem}{Theorem}[section]

\newtheorem{remark}[theorem]{Remark}
\newtheorem{proposition}[theorem]{Proposition}

\newtheorem{corollary}[theorem]{Corollary}
\newtheorem{example}[theorem]{Example}

\numberwithin{equation}{section}
\allowdisplaybreaks[1]

\theoremstyle{definition}

\newtheorem{definition}[theorem]{Definition}

\theoremstyle{remark}

\newcommand{\bA}{{\mathbf A}}
\newcommand{\bU}{{\mathbf U}}

\newcommand{\blam}{{\boldsymbol \lambda}}
\newcommand{\lam}{\lambda}
\newcommand{\bzeta}{{\boldsymbol \zeta}}
\newcommand{\bom}{{\boldsymbol \omega}}
\newcommand{\bn}{{\mathbf n}}

\newcommand{\cD}{{\mathcal D}}
\newcommand{\cF}{{\mathcal F}}
\newcommand{\cH}{{\mathcal H}}
\newcommand{\cL}{{\mathcal L}}
\newcommand{\cM}{{\mathcal M}}
\newcommand{\cR}{{\mathcal R}}
\newcommand{\cS}{{\mathcal S}}

\newcommand{\cU}{{\mathcal U}}
\newcommand{\cX}{{\mathcal X}}
\newcommand{\cY}{{\mathcal Y}}
\newcommand{\cO}{{\mathcal O}}

\newcommand{\bbZ}{{\mathbb Z}}
\newcommand{\C}{{\mathbb C}}
\newcommand{\B}{{\mathbb B}}

\newcommand{\sbm}[1]{\left[\begin{smallmatrix} #1
		\end{smallmatrix}\right]}

\begin{document}

\title[de Branges-Rovnyak spaces and commutative realizations]
{Schur-class multipliers on the Arveson space:
de Branges-Rovnyak reproducing kernel spaces and commutative transfer-function
realizations}
\author[J. A. Ball]{Joseph A. Ball}
\address{Department of Mathematics,
Virginia Tech,
Blacksburg, VA 24061-0123, USA}
\email{ball@math.vt.edu}
\author[V. Bolotnikov]{Vladimir Bolotnikov}
\address{Department of Mathematics,
The College of William and Mary,
Williamsburg VA 23187-8795, USA}
\email{vladi@math.wm.edu}
\author[Q. Fang]{Quanlei Fang}
\address{Department of Mathematics,
Virginia Tech,
Blacksburg, VA 24061-0123, USA}
\email{qlfang@math.vt.edu}

\begin{abstract}
An interesting and recently much studied generalization of the
classical Schur class is the class of contractive operator-valued
multipliers $S(\blam)$ for the reproducing kernel Hilbert space $\cH(k_{d})$ on
the unit ball ${\mathbb B}^{d} \subset {\mathbb C}^{d}$, where
$k_{d}$ is the positive kernel $k_{d}(\blam, \bzeta) = 1/(1
- \langle \blam, \bzeta \rangle)$ on ${\mathbb B}^{d}$.  The
reproducing kernel space $\cH(K_{S})$ associated with the positive
kernel $K_{S}(\blam, \bzeta) = (I - S(\blam) S(\bzeta)^{*}) \cdot
k_{d}(\blam, \bzeta)$ is a natural multivariable generalization of
the classical de Branges-Rovnyak canonical model space. A special
feature appearing in the multivariable case is that the space
$\cH(K_{S})$ in general may not be invariant under the adjoints
$M_{\lam_{j}}^{*}$ of the multiplication operators $M_{\lambda_{j}}
\colon f(\blam) \mapsto \lam_{j} f(\blam)$ on $\cH(k_{d})$.  We show
that invariance of $\cH(K_{S})$ under $M_{\lam_{j}}^{*}$ for each $j
= 1, \dots, d$ is equivalent to the existence of a realization for
$S(\blam)$ of the form $S(\blam) = D + C (I - \lam_{1}A_{1} -
\cdots - \lam_{d} A_{d})^{-1}(\lam_{1}B_{1} + \cdots + \lam_{d}
B_{d})$ such that connecting operator $\bU = \sbm{ A_{1} & B_{1} \\ \vdots &
\vdots \\ A_{d} & B_{d} \\ C & D}$ has adjoint$\bU^{*}$ which is
isometric on a certain natural subspace ($\bU$ is ``weakly
coisometric'') and has the additional
property that the state operators  $A_{1}, \dots, A_{d}$ pairwise
commute; in this case one can take the state space to be the
functional-model space $\cH(K_{S})$ and the state operators $A_{1},
\dots, A_{d}$ to be given by $A_{j} = M_{\lam_{j}^{*}|_{\cH(K_{S})}}$
(a de Branges-Rovnyak functional-model realization).
We show that this special situation always occurs for the case of inner
functions $S$ (where the associated multiplication operator $M_{S}$
is a partial isometry), and that inner multipliers are characterized
by the existence of such a realization such that the state operators
$A_{1}, \dots, A_{d}$ satisfy an additional stability property.
\end{abstract}


\maketitle

\section{Introduction}  \label{S:Intro}
\setcounter{equation}{0}

A multivariable generalization of the Szeg\"o kernel
$k(\lambda,\zeta)=(1-\lambda\bar{\zeta})^{-1}$ much studied of
late is the positive kernel
\begin{equation}
k_d(\blam,\bzeta)=\frac{1}{1-\langle \blam,  \bzeta  \rangle}
\label{1.1a}
\end{equation}
on ${\mathbb B}^{d} \times {\mathbb B}^{d}$ where
$\B^d=\left\{\blam=(\lam_1,\dots, \lam_d)\in\C^d \colon  \langle \blam,
\blam\rangle<1\right\}$ is the unit ball of the  $d$-dimensional Euclidean
space $\C^d$. By
$\langle \blam,  \bzeta \rangle=\sum_{j=1}^d \lam_j \overline{\zeta}_j$
we mean the standard inner product in $\C^d$.
The reproducing kernel Hilbert space (RKHS) $\cH(k_d)$ associated
with
$k_d$ via Aronszajn's construction \cite{aron} is a natural multivariable
analogue of the Hardy space $H^2$ of the unit disk and coincides with
$H^2$ if $d=1$.

For $\cY$ an auxiliary Hilbert space, we consider the tensor product
Hilbert space $\cH_\cY(k_d):=\cH(k_d)\otimes\cY$ whose
elements can be viewed as  $\cY$-valued functions in $\cH(k_d)$.
Then $\cH_\cY(k_d)$ can be characterized as follows:
\begin{equation}
\cH_{\cY}(k_d)=\left\{f(\blam)=\sum_{\bn \in{\mathbb Z}^d_{+}}f_{\bn}
\blam^\bn:\|f\|^{2}=\sum_{\bn \in {\mathbb Z}^d_{+}}
\frac{\bn!}{|\bn|!}\cdot \|f_{\bn}\|_{\cY}^2<\infty\right\}.
\label{char}
\end{equation}
Here and in what follows, we use standard multivariable notations: for
multi-integers $\bn =(n_{1},\ldots,n_{d})\in\bbZ_+^d$ and points
$\blam=(\lam_1,\ldots,\lam_d)\in\C^d$ we set
\begin{equation}
|\bn| = n_{1}+n_{2}+\ldots +n_{d},\qquad
\bn!  = n_{1}!n_{2}!\ldots n_{d}!, \qquad
\blam^\bn = \lambda_{1}^{n_{1}}\lambda_{2}^{n_{2}}\ldots
\lambda_{d}^{n_{d}}.
\label{mnot}
\end{equation}
By $\cL(\cU, \cY)$ we denote the space of all bounded linear operators
between Hilbert spaces $\cU$ and $\cY$. The space of multipliers
${\mathcal M}_d(\cU,\cY)$ is defined  as the space of all
$\cL(\cU,\cY)$-valued analytic functions  $S$ on $\B^d$ such that the
induced multiplication operator
\begin{equation}
M_S: \; f(\blam)\to S(\blam)\cdot f(\blam)
\label{ms}
\end{equation}
maps $\cH_{\cU}(k_d)$ into $\cH_{\cY}(k_d)$. It follows by the closed
graph theorem that for every $S\in {\mathcal M}_d(\cU,\cY)$, the
operator $M_S$ is  bounded. We shall pay particular attention to the unit
ball of  $ {\mathcal M}_{d}(\cU, \cY)$, denoted by
$$
{\mathcal S}_{d}(\cU, \cY) = \{ S \in {\mathcal M}_{d}(\cU,\cY) \colon
\| M_{S} \|_{\text{op}} \le 1 \}.
$$
Since ${\mathcal S}_{1}(\cU, \cY)$ collapses to the classical Schur class
(of holomorphic, contractive $\cL(\cU, \cY)$-valued functions on ${\mathbb
D}$),  we refer to ${\mathcal S}_{d}(\cU, \cY)$ as a generalized
($d$-variable) {\em Schur class}. Characterizations of ${\mathcal
S}_{d}(\cU, \cY)$ in terms of realizations originate to \cite{BTV, AM}.
We recall this result in the form presented in \cite{BBF2a}.
\begin{theorem}
\label{T:BTV}
Let $S$ be an $\cL(\cU, \cY)$-valued function defined
on ${\mathbb B}^{d}$.  The following are equivalent:
\begin{enumerate}
\item $S$ belongs to $\cS_d(\cU, \cY)$.
\item The kernel
\begin{equation}  \label{KS}
         K_{S}(\blam, \bzeta) = \frac{I_{\cY} - S(\blam) S(\bzeta)^{*}}{1 -
         \langle \blam, \bzeta \rangle}
\end{equation}
is positive on $\B^{d}\times\B^{d}$, i.e., there exists an operator-valued
function $H \colon{\mathbb B}^d \to \cL(\cH,\cY)$
for some auxiliary Hilbert  space $\cH$ so that
\begin{equation}
K_{S}(\blam, \bzeta) = H(\blam) H(\bzeta)^{*}.
\label{fact}
\end{equation}
\item There exists a Hilbert space $\cX$ and a
unitary connecting operator (or colligation) $\bU$ of the form
\begin{equation}
\label{1.7a}
\bU = \begin{bmatrix} A  & B \\ C & D \end{bmatrix} =
\begin{bmatrix} A_{1} & B_{1} \\ \vdots & \vdots \\ A_{d} & B_{d}
\\ C & D \end{bmatrix} \colon \begin{bmatrix} \cX \\ \cU
\end{bmatrix} \to \begin{bmatrix}\cX^{d} \\ \cY \end{bmatrix}
\end{equation}
so that $S(\blam)$ can be realized in the form
\begin{eqnarray}
S(\blam)&=&D+C\left(I_{\cX}-\lam_1A_1-\cdots-\lam_dA_d\right)^{-1}
\left(\lambda_1B_1+\ldots+\lambda_dB_d\right)\nonumber\\
&=& D + C (I - Z(\blam) A)^{-1} Z(\blam) B
\label{1.5a}
\end{eqnarray}
where we set
\begin{equation}
Z(\blam)=\begin{bmatrix}\lam_1 I_{\cX} & \ldots &
\lam_dI_{\cX}\end{bmatrix},\quad
A=\begin{bmatrix} A_1 \\ \vdots \\ A_d\end{bmatrix},\quad
B=\begin{bmatrix} B_1 \\ \vdots \\ B_d\end{bmatrix}.
\label{1.6a}
\end{equation}
\item There exists a Hilbert space $\cX$ and a
contractive connecting operator $\bU$ of the form \eqref{1.7a}
so that $S(\blam)$ can be realized in the form \eqref{1.5a}.
\end{enumerate}
\end{theorem}
In analogy with the univariate case, a realization of the form
(\ref{1.5a}) is called {\em coisometric}, {\em isometric}, {\em unitary}
or {\em contractive} if the operator $\bU$ is respectively, coisometric,
isometric, unitary or just contractive. It turns out that a more useful
analogue of ``coisometric realization'' appearing in the classical
univariate case is not that the whole connecting operator ${\bf U}^{*}$ be
isometric, but rather that ${\bf U}^{*}$ be isometric on a certain
subspace of $\cX^d\oplus\cY$.
\begin{definition}
A realization \eqref{1.5a} of $S\in\cS_d(\cU,\cY)$ is called
{\em weakly coisometric} if the adjoint $\bU^*: \, \cX^d\oplus\cY\to
\cX\oplus\cU$ of the connecting operator
is contractive and isometric on the subspace
$\begin{bmatrix}{\mathcal D}\\
\cY\end{bmatrix}\subset \begin{bmatrix}\cX^d\\ \cY\end{bmatrix}$ where
\begin{equation}  \label{domV0}
         \cD:= \operatorname{\overline{span}}
         \{Z(\bzeta)^{*}(I_{\cX} - A^{*} Z(\bzeta)^{*})^{-1} C^{*}y \colon \;
\; \bzeta\in {\mathbb B}^{d}, \, y \in \cY \} \subset \cX^{d}.
         \end{equation}
\label{D:1.2}
\end{definition}

Weakly coisometric realizations for an $S\in\cS_d(\cU,\cY)$ can be
constructed in certain canonical way as follows. Upon applying
Aronszajn's  construction to the kernel $K_S$ defined as in \eqref{KS},
(which is positive on $\B^d$ by Theorem \ref{T:BTV}), one gets the  de
Branges-Rovnyak space $\cH(K_{S})$.
A weakly coisometric
realization for $S$ with the state space equal to $\cH(K_S)$ (and output
operator $C$ equal to evaluation at zero on $\cH(K_{S})$) will be called
a {\em generalized functional-model realization}.
Here we use the term {\em generalized} functional-model realization
since it may be the case that the state space $\cH(K_{S})$ in not even
invariant under the adjoints $M_{\lam_{1}}^{*}, \dots,
M_{\lam_{d}}^{*}$ of the multiplication operators
$M_{\lam_{j}} \colon f(\blam) \mapsto \lam_{j} \cdot f(\blam)$ ($j =
1, \dots, d$) on $\cH_{\cY}(k_{d})$ and hence one cannot take the
state operators $A_{1}, \dots, A_{d}$ to be given by $A_{j} =
M_{\lam_{j}}^{*}$ as one would expect from the classical case.
As it was shown in
\cite{BBF2a}, any function $S\in\cS_d(\cU,\cY)$ admits a generalized
functional-model realization. In the univariate case, this
collapses to the well known de Branges-Rovnyak functional-model realization
\cite{dbr1, dbr2}. Another parallel to the univariate case is that
{\em any} observable  weakly coisometric realization of a Schur-class function
$S\in\cS_d(\cU,\cY)$ is unitarily equivalent to some generalized
functional-model realization (observability is a minimality condition
that is fulfilled automatically for every generalized functional-model
realization).  However, in contrast to the univariate case,
this realization is not unique in general (even up to unitary
equivalence); moreover, a function $S\in\cS_d(\cU,\cY)$ may admit
generalized functional-model realizations with the same state space
operators $A_1,\ldots,A_d$ and different input operators $B_j$'s.
A curious fact is that none of the generalized functional-model
realizations for $S$ may be coisometric.

In this paper we study another issue not present in the univariate classical
case, namely the distinction between {\em commutative realizations}
(where the state space operators $A_{1}, \dots, A_{d}$ in \eqref{1.5a}
commute with each other) versus general realizations. Commutative
realization is a natural object that appears for example in model theory
for commuting row contractions \cite{BES}: the characteristic function of
a commuting row contraction $(T_1,\ldots,T_d)$ is, by  definition, a
Schur-class function
that admits a unitary commutative
realization with the state space operators $T_1,\ldots,T_d$. It turns
out that not every $S \in \cS_d(\cU,\cY)$ can be identified as a characteristic
function of a commutative row contraction; thus not every
$S\in\cS_d(\cU,\cY)$ admits a commutative unitary
realization. Some more delicate arguments based on backward-shift
invariance in $\cH_{\cY}(k_{d})$ show that not every  $S\in\cS_d(\cU,\cY)$
admits a commutative weakly coisometric realization (see Theorem
\ref{T:comreal} below);  more surprisingly, there are Schur-class
functions that
do not admit  even {\em contractive} commutative realizations (see Example
\ref{E:3.1} below).  If the Schur-class function admits a commutative
weakly coisometric realization, then the associated de
Branges-Rovnyak space $\cH(K_{S})$ is invariant for the backward
shift operators $M_{\lam_{j}}^{*}$ and one can arrange for a
generalized functional-model realization with the additional property
that the state operators $A_{1}, \dots, A_{d}$ are given by $A_{j} =
M_{\lam_{j}}^{*}|_{\cH(K_{S})}$ for $j = 1, \dots, d$;
we say that
such a realization is a (non-generalized) {\em functional-model realization}.
The operators
$B_1,\ldots,B_d$ are not defined uniquely by $S(\blam)$,
$\bA=(A_1,\ldots,A_d)$ and $C$
(this is yet another distinction from the univariate case); however the
nonuniqueness can be described in an explicit way. Furthermore, any
observable,  commutative, weakly coisometric realization for a given $S$
is unitarily equivalent to exactly one functional-model realization
(Theorem \ref{T:canonical-universal}).

Inner functions, i.e, a Schur-class multiplier $S \in
{\mathcal  S}_{d}(\cU, \cY)$ for which the associated multiplication
operator is a partial isometry,
are special in that an inner
function necessarily has a commutative weakly
coisometric realization (see Theorem \ref{T:comreal} below).
Inner functions also play a special role as representers for
(forward) shift-invariant subspaces of $\cH_{\cY}(k_{d})$; for the
case $d=1$ this is the classical Beurling-Lax-Halmos theorem
(\cite{beurling, halmos, lax}) while the case for general $d$ appears
more recently in the work of Arveson \cite{arv} and of
McCullough-Trent \cite{mt} (for the general framework of a complete
Nevanlinna-Pick kernel).
Here we use our realization-theoretic characterization of inner multipliers
to present a new proof of the $\cH_{\cY}(k_{d})$-Beurling-Lax
theorem.  The idea in this approach is to represent the
shift-invariant subspace $\cM$ as the set of all $\cH_{\cY}(k_{d})$-solutions
of fairly general set of homogeneous interpolation conditions, and
then to construct a realization $\bU = \sbm{A & B \\ C & D}$ for
$S(\blam)$ from the operators defining the homogeneous interpolation
conditions.  For the case $d=1$, this approach can be found in
\cite{BGR} for the rational case and in \cite{BallRaney} for the
non-rational case, done there in the more complicated context
where the shift-invariant subspace $\cM$ is merely contained in the
$\cY$-valued $L^{2}$ space over the unit circle ${\mathbb T}$ and is not
necessarily contained in the Hardy space $\cH_{\cY}(k_{1}) =
H^{2}_{\cY}$.  We also use our analysis of the nonuniqueness of the
input operator $B$ in weakly coisometric realizations to characterize
the nonuniqueness in the choice of inner-function representer $S$ for
a given shift-invariant subspace $\cM$ (see Theorem \ref{T:repr}).

A more general version of the $\cH_{\cY}(k_{d})$-Beurling-Lax
theorem, where the subspace
$\cM$ is only contractively included in $\cH_{\cY}(k_{d})$ and the
representer is a not necessarily an inner Schur-class multiplier,
appears in the work of de Branges-Rovnyak \cite{dbr1, dbr2} for the
case $d=1$ and of the authors \cite{BBF1} for the case of general $d$.
The realization produced by our approach here (working with
$\cM^{\perp}$ rather than directly with $\cM$) is more explicit for
the situation where $\cM$ is presented as the solution set for a
homogeneous interpolation problem.

The paper is organized as follows.  After the present Introduction,
Section 2 recalls needed preliminaries from our earlier papers
\cite{BBF1, BBF2a} concerning weakly coisometric realizations (see
Definition \ref{D:1.2} above).  Section 3 collects the results
concerning such realizations where the collection of state operators
$A_{1}, \dots, A_{d}$ is commutative.  Section 4 specializes the
general theory to the case of inner functions.
The final Section 5 discusses connections with characteristic
functions and operator-model theory for commutative row contractions,
a topic of recent work of Bhattacharyya, Eschmeier, Sarkar and
Popescu \cite{BES, BES2, BS, Popescu-var1, Popescu-var2}, where some
extensions to more general settings are also addressed.

\section{Weakly coisometric realizations}
\label{S:coisom}

Weakly coisometric realizations of Schur-class functions are closely
related to range spaces of observability operators appearing in the
context of Fornasini-Mar\-chesini-type linear systems with evolution along
the integer lattice ${\mathbb Z}^{d}$ . Let $\bA=(A_1,\ldots,A_d)$ be a
$d$-tuple of operators in $\cL(\cX)$. If $C \in \cL(\cX, \cY)$, then the
pair $(C,{\mathbf A})$ is said
         to be an {\em output pair}.  Such an output pair is said to be {\em
         contractive} if
\begin{equation}
A_{1}^{*}A_{1} + \cdots + A_{d}^{*} A_{d} + C^{*}C\le I_{\cX},
\label{contr}
\end{equation}
to be {\em isometric} if equality holds in the above relation, and
to be {\em output-stable} if the associated observability
         operator
         \begin{equation}  \label{cO}
           \cO_{C, \bA} \colon x \mapsto C (I_{\cX} - Z(\blam) A)^{-1}x=
C (I - \lam_{1} A_{1} - \cdots -\lam_{d} A_{d})^{-1} x
          \end{equation}
(where $Z(\blam)$ and $A$ are defined as in \eqref{1.6a})
maps $\cX$ into $\cH_{\cY}(k_{d})$. As it was shown in \cite{BBF1}, any
contractive pair $(C,\bA)$ is output stable and moreover, the corresponding
observability operator ${\mathcal O}_{C, \bA}: \, \cX\to\cH_{\cY}(k_d)$ is
a contraction. An output stable pair $(C,\bA)$ is called {\em observable}
if the observability operator ${\mathcal O}_{C, \bA}$ is injective, i.e.,
$$
C (I_{\cX} - Z(\blam) A)^{-1}x\equiv 0\quad\Longrightarrow \quad x=0.
$$
Given an  output stable pair $(C,\bA)$, the kernel
\begin{equation}  \label{KCA}
K_{C,{\mathbf A}}(\blam, \bzeta):=
C (I_{\cX} - Z(\blam) A)^{-1}(I_{\cX} -A^{*}Z(\bzeta)^{*})^{-1} C^{*}.
\end{equation}
is positive on ${\mathbb B}^d\times{\mathbb B}^d$;
let $\cH(K_{C,\bA})$ denote the associated RKHS. We recall (see \cite{aron})
that any positive kernel $(\lam, \zeta) \mapsto K(\lam, \zeta) \in
\cL(\cY)$ on a set $\Omega \times \Omega$ (so $\lam, \zeta \in \Omega$)
gives rise to  a RKHS $\cH(K)$ consisting of $\cY$-valued functions on
$\Omega$ with the defining property: for each $\zeta \in \Omega$ and $y
\in \cY$, the $\cY$-valued function $K_{\zeta} y(\lambda):=K(\lam, \zeta)y$
is in $\cH(K)$ and has the reproducing property
           $$
            \langle f, K_{\zeta}y \rangle_{\cH(K)} = \langle f(\zeta), y
            \rangle_{\cY}\quad \text{for all}\quad y \in \cY, \; f\in\cH(K).
           $$
The following result from \cite{BBF1} gives the close connection
between spaces of the form $\cH(K_{C, \bA})$ and ranges of
observability operators.
\begin{theorem}  \label{T:3-1.2nc} (See \cite[Theorem 3.20]{BBF1}.)
Let $(C,\bA)$ be a contractive pair with $C\in{\mathcal L}(\cX,\cY)$
and with associated positive kernel $K_{C, \bA}$ given by \eqref{KCA} and
the observability operator $\cO_{C, \bA}$ given by \eqref{cO}.  Then:
         \begin{enumerate}
\item The reproducing kernel Hilbert space $\cH(K_{C,\bA})$ is
characterized as
$$
\cH(K_{C,\bA}) = \operatorname{Ran}\, \cO_{C, \bA}
$$
with the lifted norm given by
\begin{equation}
     \| \cO_{C, \bA} x \|_{\cH(K_{C, \bA})} = \| Q x \|_{\cX}
\label{june26}
\end{equation}
where $Q$ is the orthogonal projection onto
$(\operatorname{Ker}\, \cO_{C, \bA})^{\perp}$.
        \item The operator ${\cO}_{C,\bA}$
               is a contraction of $\cX$ into $\cH(K_{C,\bA})$.
               It is an isometry if and only if the the pair $(C,\bA)$ is
               observable.
\item The space $\cH(K_{C,\bA})$ is contractively included in the
Arveson space $\cH_{\cY}(k_d)$; it is isometrically included
in $\cH_{\cY}(k_d)$ if and only if ${\cO}_{C,\bA}$
(as an operator from $\cX$ into $\cH_{\cY}(k_d)$) is a partial isometry.
               \end{enumerate}
              \end{theorem}
If $S$ is realized as in \eqref{1.5a} and  $\bU$ is the connecting
operator given by \eqref{1.7a}, then the associated kernels $K_S$ and
$K_{C, \bA}$ (defined in \eqref{KS} and \eqref{KCA}, respectively)
are related by the following easily verified identity
\begin{align}
K_S(\blam, \bzeta)&=:K_{C,{\mathbf A}}(\blam, \bzeta)\label{1.6u}\\
&+\begin{bmatrix}C (I-Z(\blam)A)^{-1}Z(\blam) & I\end{bmatrix}
\frac{I-\bU\bU^*}{1-\langle \blam, \bzeta \rangle}
\begin{bmatrix}Z(\bzeta)^*(I_{\cX} -
A^{*}Z(\bzeta)^{*})^{-1} C^{*}\\ I\end{bmatrix}\notag
\end{align}
and then, it is easily shown (see Proposition 1.5 in \cite{BBF2a}
for details) that the second term on the right vanishes if and only if
${\bf U}^{*}$ is isometric on the space ${\mathcal D}\oplus\cY$ defined as
in Definition \ref{D:1.2}. This observation leads us to the following
intrinsic kernel characterization as to when a given contractive
realization is a weakly coisometric realization.

\begin{proposition} \label{R:1.3}
A contractive realization \eqref{1.5a} of $S\in\cS_d(\cU,\cY)$ is  weakly
coisometric if and only if the kernel  $K_{S}(\blam, \bzeta)$ associated
to $S$ via \eqref{KS}  can alternatively be written as
\begin{equation}
K_{S}(\blam, \bzeta) = K_{C, \bA}(\blam, \bzeta)
\label{2.1}
\end{equation}
where $K_{C, \bA}$ is given by \eqref{KCA}.
\end{proposition}

Proposition \ref{R:1.3} states that once a contractive realization
$\bU= \sbm{A & B \\ C & D}$ of $S$ is such that \eqref{2.1} holds, then
this realization is weakly coisometric. The next result asserts that
equality \eqref{2.1} itself guarantees the existence of weakly
coisometric realizations for $S$ with preassigned $C$ and
$\bA=(A_1,\ldots,A_d)$.

\begin{theorem}  \label{T:CAtoS} (See \cite[Theorem 2.4]{BBF2a}.)
Suppose that a Schur-class function $S \in \cS_{d}(\cU, \cY)$ and
a contractive pair $(C, \bA)$ are such that \eqref{2.1} holds and let
$D:=S(0)$. Then there exist operators $B \colon \cU \to \cX^d$ so
that the operator $\bU$ of the form \eqref{1.7a} is
weakly coisometric and $S$ can be realized as in \eqref{1.5a}.
\end{theorem}

The pair $(C, \bA)$  for a weakly coisometric realization can be
constructed in a certain canonical way.

\begin{theorem} \label{T:noncomreal} (See \cite[Theorem 3.20]{BBF2a}.)
Let $S \in \cS_{d}(\cU, \cY)$ and let $\cH(K_{S})$ be the associated de
Branges-Rovnyak space. Then:
\begin{enumerate}
\item There exist bounded operators $A_{j} \colon \cH(K_{S})
\to \cH(K_{S})$ such that
\begin{equation}
f(\blam) - f(0)  = \sum_{j=1}^{d} \lam_{j} (A_{j}f)(\blam)
\quad\mbox{for every $f \in \cH(K_{S})$ and $\blam \in {\mathbb B}^{d}$},
\label{dop1}
\end{equation}
and
\begin{equation}
        \sum_{j=1}^{d} \|A_{j}f\|^{2}_{\cH(K_{S})}  \le
             \|f\|^{2}_{\cH(K_{S})} - \| f(0)\|^{2}_{\cY}.
\label{dop2}
\end{equation}
\item  There is a weakly coisometric realization \eqref{1.5a} for $S$
with state space $\cX$ equal to $\cH(K_{S})$ with the state operators
$A_1,\ldots, A_d$ from part (1) and the operator $C \colon  \cH(K_S)\to \cY$
defined by
\begin{equation}
Cf=f(0)  \quad\text{for all}\quad f\in\cH(K_S).
\label{dop8}
\end{equation}
\end{enumerate}
\end{theorem}
Equality \eqref{dop1} means that the operator tuple
${\bf A}=(A_1,\ldots,A_d)$ solves {\em the Gleason problem} \cite{gleason}
for $\cH(K_{S})$. Let us say that ${\bf A}$ {\em is a
contractive solution of the Gleason problem} if in addition relation
\eqref{dop2} holds for every $f\in\cH(K_S)$ or, equivalently, if the pair
$(C,{\bf A})$ is contractive where $C: \, \cH(K_S)\to\cY$ is defined as
in \eqref{dop8}. Theorem \ref{T:noncomreal} shows that any contractive
solution ${\mathbf A} = (A_{1},\dots, A_{d})$ of the Gleason problem for
$\cH(K_{S})$ gives rise to a weakly coisometric realization for
$S \in\cS_{d}(\cU, \cY)$ (not unique, in general). Let us call any such
weakly coisometric realization  a {\em generalized functional-model}
realization of $S(\blam)$.  We note that any generalized functional-model
realization of $S$ is observable and that the formula
   \begin{equation}  \label{H-canonical}
         K_{S}(\cdot, \bzeta) y = (I - A^{*} Z(\bzeta)^{*})^{-1} C^{*} y
\quad (y\in\cY, \; \zeta\in{\mathbb B}^d)
\end{equation}
is valid for any generalized functional-model realization. Furthermore,
if
\begin{equation}
{\mathbf U} = \begin{bmatrix} A & B \\ C & D \end{bmatrix} \colon
\begin{bmatrix} \cH(K_{S}) \\ \cU \end{bmatrix} \to \begin{bmatrix}
          \cH(K_{S})^{d} \\ \cY \end{bmatrix}
\label{gfm}
\end{equation}
is a generalized functional model realization for an $S \in {\mathcal
S}_{d}(\cU, \cY)$ then the space ${\mathcal D}$
introduced in \eqref{domV0} can be described in the following explicit
functional form
\begin{equation}
{\mathcal D}=\operatorname{\overline{span}} \{
              Z(\bzeta)^{*} K_{S}(\cdot, \bzeta) y \colon \; \bzeta \in
{\mathbb
              B}^{d},\, y \in \cY \}.
\label{canonical1}
\end{equation}
Then a simple calculation shows that ${\mathcal D}^{\perp}=\cH(K_{S})^{d}
\ominus {\mathcal D}$ can be characterized in similar terms as
\begin{equation}
{\mathcal D}^{\perp}  = \{ h  \in \cH(K_{S})^{d} \colon \;
             Z(\blam) h(\blam) \equiv 0 \}.
             \label{canonical2}
\end{equation}

\section{Realizations with commutative state space operators}
\label{S:com}

Schur-class functions that admit unitary realizations of the form \eqref{1.5a}
with commutative state space tuple $\bA=(A_{1}, \dots, A_{d})$ is a
natural object appearing in the model theory for
commutative row contractions (see \cite{BES}): the characteristic function
of a commutative row contraction (see formula \eqref{charf} below) is
a Schur-class function of this type (subject to an additional normalization).
In the commutative context, a key role is played by the commuting
$d$-tuple ${\mathbf  M}_{\blam}:=(M_{\lam_1},\ldots,M_{\lam_d})$
consisting of operators of multiplication by the coordinate functions of
$\C^d$ which will be called {\em the shift}
(operator-tuple) of $\cH_{\cY}(k_d)$, whereas the commuting $d$-tuple
${\bf M}^*_\blam:=(M^*_{\lam_1},\ldots,M^*_{\lam_d})$ consisting
of the adjoints of $M_{\lam_j}$'s (in the metric of $\cH_{\cY}(k_d)$)
will be referred to as to the {\em backward shift}. By the
characterization \eqref{char} and in notation \eqref{mnot}, the monomials
$\frac{\bn!}{|\bn|!}\blam^\bn$ form an orthonormal basis in $\cH(k_d)$
and then a simple calculation shows that
\begin{equation}
M_{\lambda_j}^* \blam^{\bf m}=\frac{m_j}{|{\bf m}|}\blam^{{\bf m}-e_j} \;
\;
(m_j\ge 1)\quad\mbox{and}\quad M_{\lambda_j}^* \blam^{\bf m}=0 \; \;
(m_j=0)
\label{bs}
\end{equation}
where ${\bf m}=(m_1,\ldots,m_d)$ and $e_j$ is the $j$-th
standard coordinate vector of $\C^d$. Some properties of the
shift tuple ${\bf M}^*_\blam$ needed in the sequel are listed below
(for the proof, see e.g., \cite[Proposition 3.12]{BBF1}). In the
formulation and in what follows, we use multivariable power notation
$$
\bA^\bn:= A_{1}^{n_{1}}A_{2}^{n_{2}}\ldots A_{d}^{n_{d}}
$$
for any $d$-tuple ${\mathbf A} = (A_{1}, \dots, A_{d})$ of commuting
operators on a space $\cX$ and any $\bn=(n_1,\ldots,n_d)\in{\mathbb
Z}_+^d$.
\begin{proposition}  \label{P:5.1}
Let ${\bf M}_{\blam}^*$ be the $d$-tuple of backward shifts on
$\cH_{\cY}(k_d)$ and let
\begin{equation}  \label{G}
G: \; f\mapsto f(0)\quad(f \in \cH_{\cY}(k_{d}))
\end{equation}
be the operator of evaluation at the origin. Then:
\begin{enumerate}
             \item
             For every $f \in \cH_{\cY}(k_{d})$ and every $\blam \in
             {\mathbb B}^{d}$ we have
             \begin{equation}
             f(\blam)-f(0)=\sum_{j=1}^d \lambda_j(M_{\lambda_j}^*f)(\blam).
             \label{5.12}
             \end{equation}

             \item The pair $(G,{\bf M}_{\blam}^*)$ is isometric and
the  associated observability operator is the identity operator:
\begin{equation}  \label{G=I}
{\cO}_{G, {\bf M}_{\blam}^*} =I_{\cH_{\cY}(k_d)}.
\end{equation}
\item
The $d$-tuple ${\bf M}_{\blam}^*$ is strongly stable,
that is,
\begin{equation}
\label{5.13a}
\lim_{N \to \infty} \sum_{\bn \in\bbZ^{d}_{+}
\colon |\bn| = N} \frac{N!}{\bn!} \|({\bf M}_{\blam}^*)^\bn f
\|_{\cH_{\cY}(k_d)}^{2} =0\quad\mbox{for every}\quad f \in\cH_{\cY}(k_d).
\end{equation}
\end{enumerate}
\end{proposition}
We will also need the commutative analogue of Theorem \ref{T:3-1.2nc}
(see \cite[Theorem 3.15]{BBF1} for the proof).

             \begin{theorem} \label{T:BBF1}
             Let $(C, \bA)$ be a contractive pair such that $C\in{\mathcal
L}(\cX,\cY)$ and the $d$-tuple $\bA=(A_1,\ldots,A_d)\in{\mathcal L}(\cX)^d$
is commutative. Let  $K_{C, \bA}$ be the associated kernel given by
\eqref{KCA}. Then:
\begin{enumerate}
\item The reproducing kernel Hilbert space $\cH(K_{C,\bA})$
is invariant under $M_{\lam_j}^*$ for $j=1,\ldots,d$ (${\bf
M}_\blam^*$-invariant) and  the difference-quotient inequality
$$
                \sum_{j=1}^d \|M_{\lambda_j}^*f\|^2_{\cH(K_{C,\bA})}
                \le\|f\|^2_{\cH(K_{C,\bA})}-\|f(0)\|^2_{\cY}
$$
                holds for every $f\in\cH(K_{C,\bA})$.
\item The space $\cH(K_{C,\bA})$ is contractively included in
$\cH_{\cY}(k_{d})$. The inclusion is isometric exactly when the pair $(C,
\bA)$ is isometric:
\begin{equation}  \label{4.15}
          I_{\cX} - A_{1}^{*}A_{1} - \cdots - A_{d}^{*}A_{d} = C^{*}C,
\end{equation}
and $\bA$ is strongly stable:
         \begin{equation}  \label{stable}
            \lim_{N \to \infty} \sum_{\bn \in {\mathbb Z}^{d}_{+} \colon
            |\bn| = N} \frac{N}{\bn !} \| \bA^{\bn} x \|_{\cX}^{2} = 0
\quad\text{for all}\quad x \in \cX.
            \end{equation}
\end{enumerate}
\end{theorem}
If one drops the requirement of the connecting operator $\bU$ being
contractive,  constructing a commutative realization is not an issue
not only for Schur-class functions, but even for functions from
$\cH(k_d)\otimes \cL(\cU,\cY)$. Indeed, for an $S\in\cH(k_d)\otimes
\cL(\cU,\cY)$, let
$$
C=G,\quad D=S(0),\quad
A_j=M_{\lambda_j}^*,\quad B_j=M_{\lambda_j}^*M_S\vert_{\cU}
\quad (j=1,\ldots,d)
$$
where $M_S\vert_{\cU}: \, u\mapsto S(\blam)u$ and $G\colon
\cH_{\cY}(k_d)\to\cY$ is given by  \eqref{G}. Pick a vector $u\in\cU$
and note that on account of  equality \eqref{G=I} and equality
\eqref{5.12} applied to $f(\blam)=S(\blam)u$,
\begin{eqnarray*}
(D + C(I - Z(\blam)A)^{-1}Z(\blam)B)u&=&
S(0)u+{\cO}_{G, {\bf
M}_{\blam}^*}(Z(\blam)Bu)\\
&=&S(0)u+Z(\blam)Bu\\ &=&
S(0)u+\sum_{j=1}^d \lambda_j(M_{\lambda_j}^*Su)(\blam)\\
&=&S(0)u+S(\blam)u-S(0)u=S(\blam)u
\end{eqnarray*}
and thus, ${\bf U}=\sbm {A & B \\ C & D}$ is a realization for $S$.
This realization is commutative and observable. However, it is not
contractive: a simple calculation
based again on identity \eqref{5.12} shows that for $\bU$ as above
and $g=\begin{bmatrix}f \\ u\end{bmatrix}
\in\begin{bmatrix} \cH_{\cY}(k_d) \\ \cU\end{bmatrix}$, we have
$$
\|g\|^2-\|\bU
g\|^2=\|f\|_{\cH_{\cY}(k_d)}^2+\|u\|^2_\cU-\|f+Su\|_{\cH_{\cY}(k_d)}^2
$$
which cannot be nonnegative for all $f\in\cH_{\cY}(k_d)$ and $u\in\cU$
unless $S(\blam)\equiv 0$.

Since our primary object of interest are  Schur-class functions for which
norm-constrained (contractive, unitary and all intermediate) realizations
do exist (by Theorem \ref{T:BTV}), it is natural to construct
commutative realizations of the same types. Note that Theorem \ref{T:BTV}
and more specific Theorem \ref{T:noncomreal} give no clue as to when
and how one can achieve a such a realization of a given $S\in\cS_d(\cU,\cY)$.
The next proposition shows that there are Schur-class functions which do not
have a commutative contractive realization.
\begin{proposition}
\label{P:doesnot}
Let $S\in\cS_d(\cU,\cU)$ be such that the associated de Branges-Rovn\-yak
space $\cH(K_S)$ is finite-dimensional and is not ${\bf
M}_\blam^*$-invariant. Then $S$ does not have a commutative contractive
realization.
\end{proposition}
\begin{proof}
Assume that $S$ admits a contractive commutative realization
\eqref{1.5a} with a contractive ${\bf U}=\sbm{ A&B\\ C&D}$. Since
${\bf U}$ is contractive, the formula
\eqref{1.6u} for $K_S$ can be written in the form
            $$
K_S(\blam, \bzeta)=
K_{C,\bA}(\blam,\bzeta)+
            \frac{S_1(\blam)S_1(\bzeta)^*}{1-\langle \blam, \bzeta \rangle}
            $$
            where $S_1\in\cS_d(\cF,\cU)$ is a Schur-class function with an
appropriately chosen coefficient space $\cF$ (the explicit formula
for $S_1$ is not that important). If $S_1\not\equiv 0$, then  ${\mathcal
H}(K_S)$ contains $S_1{\mathcal H}_{\cF}(k_d)$ and therefore is infinite
dimensional which contradicts one of the assumptions about $\cH(K_S)$. If
$S_1\equiv 0$, then $K_S(\blam, \bzeta)=K_{C, \bA}(\blam, \bzeta)$
and therefore $\cH(K_S)=\cH(K_{C,\bA})$. Since the tuple $\bA$ is
commutative, the space $\cH(K_S)=\cH(K_{C,\bA})$ is
${\bf M}_\blam^*$-invariant (by Theorem \ref{T:BBF1}) which contradicts
another assumption about $\cH(K_S)$.
\end{proof}

\begin{example}  \label{E:3.1}
{\rm For a concrete example of a Schur-class function satisfying assumptions
in Proposition \ref{P:doesnot}, let
\begin{equation}  \label{ex3.1}
S(\lam_1,\lam_2)=\frac{1}{4-\lam_1 \lam_2}\left[\begin{array}{cccc}
          2\sqrt{3} \lam_1 & \sqrt{3}\lam _2^2 & 2-2\lam_1 \lam_2 & -3 \lam_2 \\
          \sqrt{3} \lam_1^2 & 2\sqrt{3} \lam_2 & -3 \lam_{1} & 2-2 \lam_1
          \lam_2
          \end{array}\right].
          \end{equation}
A straightforward calculation gives
          \begin{eqnarray}
          K_S(\blam , \bzeta)&:=&
          \frac{I_2-S(\blam)S(\bzeta)^*}{1-\lam_1\overline{\zeta}_1-
          \lam_2\overline{\zeta}_2}\nonumber\\
          &=&\frac{3}{(4-\lam_1 \lam_2)(4-\overline{\zeta}_1\overline{\zeta}_2)}
          \left[\begin{array}{cc}
          2 & \lam_2 \\ \lam_1 & 2\end{array}\right]\left[\begin{array}{cc}
          2 & \overline{\zeta}_1 \\ \overline{\zeta}_2 &
          2\end{array}\right].\nonumber
          \end{eqnarray}
          Thus the kernel $K_S(\blam, \bzeta)$ is positive on
$\B^2\times\B^2$ and $S\in{\mathcal S}_2({\mathbb C}^4,{\mathbb C}^2)$.
The associated de Branges-Rovnyak space ${\mathcal H}(K_S)$ is
spanned by rational functions
          $$
          f_1(\blam)=\frac{4}{4-\lam_1 \lam_2}\left[\begin{array}{c}2\\
\lam_1\end{array}\right]
          \quad\mbox{and}\quad
          f_2(\blam )=\frac{4}{4-\lam_1 \lam_2}\left[\begin{array}{c}\lam_2\\
          2\end{array}\right].
          $$
          Furthermore, since by \eqref{bs} we have
          $$
          M_{\lam_1}^*(\lam_1^{n_1}\lam_2^{n_2})=
          \frac{n_1}{n_1+n_2}\lam_1^{n_1-1}\lam_2^{n_2},
          $$
         it follows that
\begin{equation}  \label{notrational}
          M_{\lam_1}^*\left(\frac{4\lam_1}{4-\lam_1 \lam_2}\right)=
          M_{\lam_1}^*\left(\sum_{j=0}^\infty
          \frac{\lam_1^{j+1}\lam_2^j}{4^j}\right)=
          \sum_{j=0}^\infty\frac{j+1}{2j+1}\left(\frac{\lam_1  \lam_2}{4}
          \right)^{j}.
\end{equation}

The latter function is rational if and only if the single-variable
function $F(\lam) = \sum_{j=0}^{\infty} \frac{j+1}{2j+1} \lam^{j}$ is
rational.  By the well-known Kronecker theorem, $F$ in turn is
rational if and only if the associated infinite Hankel matrix
$$ {\mathbb H} = [s_{i+j}]_{i,j=0}^{\infty} \text{ where }
       s_{k} = \frac{k+1}{2k+1}
$$
has finite rank. However one can check that the finite Hankel matrices
${\mathbb H}_{n} = [s_{i+j}]_{i,j=0}^{n}$ have full rank for all
$n=0,1,2,\dots$ and hence $F(\lam)$ is not rational.  We conclude that
the function on the right hand side in \eqref{notrational} is not rational.
Now it follows that $M_{\lam_1}^*f_1$ does
not belong to ${\mathcal H}(K_S)$ . Therefore ${\mathcal H}(K_S)$ is
not invariant  under $M_{\lam_1}^*$ and since $\dim \, \cH(K_S)=2<\infty$, the
function  $S$ does not admit contractive
commutative realizations by Proposition \ref{P:doesnot}.}
\end{example}
A characterization  of which Schur-class functions do admit contractive
commutative realizations will be given in Theorem \ref{T:3.5a} below.
The next result gives a characterization of Schur-class functions
that admit weakly coisometric commutative realizations.

\begin{theorem}  \label{T:comreal}
A Schur-class function $S\in\cS_d(\cU,\cY)$ admits a commutative weakly
coisometric realization if and only if the following conditions hold:
\begin{enumerate}
           \item The
associated de Branges-Rovnyak space $\cH(K_S)$ is ${\bf M}_\blam^*$-invariant,
and
\item the inequality
\begin{equation}
\sum_{j=1}^d\|M_{\lam_j}^*f\|^2_{\cH(K_S)}\le
\|f\|^2_{\cH(K_S)}-\|f(0)\|^2_{\cY}\quad\mbox{holds for all}\quad f \in
\cH(K_{S}).
\label{3.1}
\end{equation}
\end{enumerate}
\end{theorem}

\begin{proof}
            To prove necessity, suppose that $S \in \cS_{d}(\cU, \cY)$ admits
            a weakly coisometric realization \eqref{1.5a}. As noted in
Proposition \ref{R:1.3}, it
            follows that $\cH(K_S) = \cH(K_{C, \bA})$.  Since $\bA$ is
            commutative, Theorem \ref{T:BBF1} implies that the space
            $\cH(K_S) = \cH(K_{C, \bA})$ is ${\mathbf
            M}_{\blam}^{*}$-invariant with the inequality \eqref{3.1} holding.

To prove sufficiency, suppose that $S \in \cS_{d}(\cU, \cY)$ is such that
$\cH(K_{S})$ is ${\mathbf M}_{\blam}^{*}$-invariant with \eqref{3.1}
holding.  Define operators $A_1,\ldots,A_d \colon \; \cH(K_{S})\to
\cH(K_{S})$, $\; C \colon \; \cH(K_{S}) \to \cY\;$ and $\;D\colon \;
\cU\to \cY$ by
\begin{equation}
A_{j} = M_{\lam_{j}}^{*}\vert_{\cH(K_{S})} \;\; (j=1,\ldots,d),\quad
C \colon \; f \to f(0),\quad D = S(0).
\label{defCA}
\end{equation}
Formula  \eqref{5.12} tells us that the operators $M_{\lam_{1}}^{*},
\ldots,M_{\lam_{d}}^{*}$ solve the Gleason problem for $\cH_{\cY}(k_d)$.
In particular, restriction of this formula to $f \in \cH(K_{S})$
can be written in terms of the operators \eqref{defCA} in the form
\eqref{dop1}, which means that $A_1,\ldots,A_d$ solve the Gleason problem
for $\cH(K_S)$. Then we apply Theorem \ref{T:noncomreal} (part (2)) to
conclude that there is a choice of $B\colon \cU \to \cH(K_S)^d$
with $\bU$ of the form \eqref{gfm}
weakly coisometric so that $S(\blam) = D + C(I - Z(\blam)A)^{-1}
Z(\blam)B$. This completes the proof.
             \end{proof}

Note that the proof of Theorem \ref{T:comreal} obtains
a realization for $S \in
          \cS_{d}(\cU, \cY)$ of a special form under the assumption that
          $\cH(K_{S})$ is ${\mathbf M}_{\lam}^{*}$-invariant:  the state
space $\cX$ is  taken to be the de Branges-Rovnyak space $\cH(K_{S})$ and
          the operators $\bA= (A_{1}, \dots, A_{d})$, $C$, $D$ are given by
          \eqref{defCA}; only the operators $B_{j} \colon \cU \to \cH(K_{S})$
          remain to be determined.  We shall say that {\em any} contractive
          realization of a given Schur-class function
          $S$ of this form (i.e., with $\cX = \cH(K_{S})$ and $A,C,D$ given
          by \eqref{defCA}) is a {\em functional-model
          realization} of $S$. It is readily seen that any functional-model
realization is also a generalized functional-model realization; in
particular, it is weakly coisometric and observable.

Let us recall that two colligations
$$
{\mathbf U} = \begin{bmatrix}A & B \\ C & D\end{bmatrix}
\colon\cX \oplus \cU \to \cX^{d} \oplus \cY\quad\mbox{and}\quad
\widetilde{{\mathbf U}}
= \begin{bmatrix}\widetilde{A}& \widetilde{B} \\ \widetilde{C} &
\widetilde{D}\end{bmatrix}
\colon\widetilde{\cX} \oplus \cU \to \widetilde{\cX}^{d}
\oplus \cY
$$
are said to be {\em unitarily equivalent} if
there is a unitary operator $U \colon \cX \to \widetilde{\cX}$ such that
$$
\begin{bmatrix} \oplus_{k=1}^{d} U & 0 \\ 0 & I_{\cY} \end{bmatrix}
\begin{bmatrix} A & B \\ C & D \end{bmatrix} =
\begin{bmatrix}\widetilde{A}&
\widetilde{B} \\ \widetilde{C} &
\widetilde{D} \end{bmatrix} \begin{bmatrix} U & 0 \\ 0
& I_{\cU} \end{bmatrix}.
$$
As it was shown in \cite{BBF2a}, any observable  weakly coisometric
realization of a Schur-class function $S\in\cS_d(\cU,\cY)$ is unitarily
equivalent to some generalized functional-model realization.
An analogous result concerning the universality of functional-model
realizations among commutative realizations is more specific.

\begin{theorem}
\label{T:canonical-universal}
Suppose that $S(\blam) \in \cS_{d}(\cU, \cY)$ is a Schur-class
function that admits functional-model realizations.  Then any commutative,
observable, weakly coisometric realization of $S$ is unitarily equivalent
to exactly one functional-model realization of $S$.
\end{theorem}
\begin{proof}
Let
$S(\blam)=D+\widetilde{C}(I_{\widetilde{\cX}}-Z(\blam)\widetilde{A})^{-1}
Z(\blam)\widetilde{B}$ be a commutative, observable, weakly coisometric
realization of $S$. Then $K_S(\blam, \bzeta)=K_{\widetilde C, \widetilde
\bA}(\blam, \bzeta)$ by Proposition \ref{R:1.3}. Define operators
$A_j$'s and $C$ as in \eqref{defCA}. Since $S$ admits functional-model
realizations (that contain $A_j$'s and $C$ and are weakly coisometric),
then we have also $K_{C, \bA}(\blam, \bzeta) =K_S(\blam, \bzeta)$.
Therefore $K_{C, \bA}=K_{\widetilde C, \widetilde \bA}$. Since the
the pairs $(\widetilde{C}, \widetilde \bA)$ and
$(C,\bA)$ are observable, the latter equality implies (see  \cite[Theorem
3.17]{BBF1}) that there exists a unitary operator
$U\colon \cH(K_S)\to\widetilde{\cX}$ such that
$$
         C=\widetilde{C}U\quad\mbox{and}\quad A_j=U^*\widetilde{A}_jU
         \quad\mbox{for} \; j=1,\ldots,d.
$$
Now we let $B_j:=U^*\widetilde{B}_j: \, \cY\to \cH(K_S)$ for
$j=1,\ldots,d$ which is the unique choice that guarantees the realization
$\bU=\sbm{A & B \\ C & D}$ (with $A$ and $B$ defined as in \eqref{1.6a})
to be  unitarily equivalent to the original realization
$\widetilde{\bU}=\sbm{\widetilde{A} & \widetilde{B} \\ \widetilde{C} &
D}$ and it is functional-model realization due to the canonical
choice of $C$ and $A_j$'s.
\end{proof}
\begin{corollary}  \label{C:uneq}
Let $\bU^\prime=\sbm{A^\prime & B^\prime \\ C^\prime & D}$ and
$\bU^{\prime\prime}=\sbm{A^{\prime\prime} & B^{\prime\prime} \\
C^{\prime\prime} &
D}$ be two observable commutative weakly coisometric realizations
         of a Schur-class function $S\in\cS_d(\cU,\cY)$. Then  the pairs
$(C^{\prime}, \bA^{\prime})$ and $(C^{\prime \prime},
\bA^{\prime \prime})$ are unitarily equivalent.
\end{corollary}

\begin{proof}
By Theorem \ref{T:canonical-universal},
the pairs $(C^\prime, \bA^\prime)$ and $(C^{\prime\prime},
\bA^{\prime\prime})$
are both unitarily equivalent to the canonical pair $(C,\bA)$ with
$C$ and $\bA=(A_1,\ldots,A_d)$ defined as in \eqref{defCA}.
Hence $(C^{\prime}, \bA^{\prime})$ and $(C^{\prime \prime},
\bA^{\prime \prime})$ are unitarily equivalent to each other.
\end{proof}

\begin{remark}
\label{R:3.7}
{\rm It was pointed out in \cite{BBF2a} and justified by examples (e.g.,
\cite[Example 3.5]{BBF2a}) a Schur-class function may have many weakly
coisometric observable realizations with associated output pairs
$(C,\bA)$  not unitarily equivalent . Theorem \ref{T:canonical-universal} above
shows that if $S\in\cS_d(\cU,\cY)$ admits a commutative weakly
coisometric realization, then the output pair $(C,\bA)$ of {\em any}
commutative weakly coisometric observable realization is uniquely
defined up to unitary equivalence. The example below shows that
in the latter case, $S$ may also admit many {\em noncommutative}
observable  weakly coisometric realizations with  output pairs
not unitarily equivalent. This example is of certain interest
because of Theorem \ref{T:doesnot} below showing that this
situation is not relevant if $S$ is an inner multiplier.}
\end{remark}
\begin{example}
\label{E:3.8}
{\rm Take the matrices
\begin{equation}\label{aug1}
C=\begin{bmatrix}\frac{1}{2} & 0 & 0 \end{bmatrix},  \quad
A_{0,1}=\begin{bmatrix} 0 & 1 & 0 \\ 0 & 0 & 0 \\ 0
&0 &0\end{bmatrix},  \quad
A_{0,2}=\begin{bmatrix} 0 & 0 & 1 \\ 0 & 0 & 0
       \\ 0 & 0 & 0\end{bmatrix},
\end{equation}
\begin{equation}\label{aug2}
B_{0,1}=\begin{bmatrix}0 & 0 & 0 & 0 & 0  \\
1 & 0 & 0 & 0 & 0 \\
0 & 1 & 0 & 0 & 0\end{bmatrix},\qquad
B_{0,2}=\begin{bmatrix}0 & 0 & 0 & 0 & 0  \\
0 & 0 & 1 & 0 & 0 \\
0 & 0 & 0 & 1 & 0\end{bmatrix},
\end{equation}
\begin{equation}
D=\begin{bmatrix}  0 & 0 & 0 & 0 &\frac{\sqrt{3}}{2} \end{bmatrix}
\label{aug4}
\end{equation}
so that the $7\times 8$ matrix
$$
{\bf U}_0=\begin{bmatrix}A_{0,1} & B_{0,1} \\ A_{0,2} & B_{0,2} \\
C & D\end{bmatrix}
$$
is coisometric. Then the characteristic function of the colligation ${\bf
U}_0$,
\begin{equation}\label{aug6}
S(\blam)
=D+C(I-\lambda_1A_{0,1}-\lambda_2A_{0,2})^{-1}
(\lambda_1B_{0,1}+\lambda_2B_{0,2})
\end{equation}
belongs to the Schur class $\cS_2(\C^5,\C)$. It is readily seen that
\begin{equation}\label{aug7}
C(I-\lambda_1A_{0,1}-\lambda_2A_{0,2})^{-1}=
\frac{1}{2}\cdot \begin{bmatrix} 1 & \lambda_1 & \lambda_2\end{bmatrix}
\end{equation}
which being substituted along with \eqref{aug2}, \eqref{aug4} into
\eqref{aug6} gives the explicit formula
\begin{equation}\label{aug8}
S(\blam)=\frac{1}{2}\cdot
\begin{bmatrix}\lambda_1^2 & \lambda_1\lambda_2 & \lambda_1\lambda_2 &
\lambda_2^2 & \sqrt{3}\end{bmatrix}.
\end{equation}
It is readily seen that the pair $(C,\bA_0)$ is observable (where we
let $\bA_{0} = (A_{0,1}, A_{0,2})$) and thus,
representation \eqref{aug6} is a coisometric (and
therefore, also weakly coisometric) observable realization of
the function $S\in\cS_2(\C^{5},\C)$ given by \eqref{aug8}. Then we also
have
\begin{eqnarray}
K_S(\blam,\bzeta)&=C(I-\lambda_1A_{0,1}-\lambda_2A_{0,2})^{-1}
(I-\bar{\zeta}_1A_{0,1}^*-\bar{\zeta}_2A_{0,2}^*)^{-1}C^*\nonumber\\
&=K_{C,\bA_0}(\blam,\bzeta).\label{aug9}
\end{eqnarray}
Now let us consider the matrices
\begin{equation}\label{aug10}
A_{\gamma,1}=\begin{bmatrix} 0 & 1 & 0 \\ 0 & 0 & 0 \\
\gamma &0 &0\end{bmatrix}\quad\mbox{and}\quad
A_{\gamma,2}=\begin{bmatrix} 0 & 0 & 1 \\ -\gamma & 0
& 0     \\ 0 & 0 & 0\end{bmatrix}
\end{equation}
where $\gamma\in\C$ is a parameter, and note that
$$
C(I-\lambda_1A_{\gamma,1}-\lambda_2A_{\gamma,2})^{-1}=
\frac{1}{2}\cdot \begin{bmatrix} 1 & \lambda_1 & \lambda_2\end{bmatrix}
$$
for every $\gamma$. In particular, the pair  $(C,\bA_\gamma)$ is
observable for every $\gamma$. The latter equality together with
\eqref{aug9} gives
\begin{equation}\label{aug11}
K_S(\blam,\bzeta)=K_{C,\bA_\gamma}(\blam,\bzeta).
\end{equation}
Now pick any $\gamma$ so that $|\gamma|<\sqrt{\frac{3}{8}}$. As it is
easily seen, the latter inequality is equivalent to the pair $(C,{\bf
A}_\gamma)$ being contractive. Thus, we have a Schur-class function $S$
and a contractive pair $(C,{\bf A}_\gamma)$ such that equality
\eqref{aug11} holds. Then by Theorem \ref{T:CAtoS}, there exist operators
$B_{\gamma,1}$ and $B_{\gamma,2}$ so that the operator
$$
{\bf U}_\gamma=\begin{bmatrix}A_{\gamma,1} & B_{\gamma,1} \\ A_{\gamma,2}
& B_{\gamma,2} \\ C & D\end{bmatrix}
$$
is weakly coisometric and $S$ can be realized as
$$
S(\blam)
=D+C(I-\lambda_1A_{\gamma,1}-\lambda_2A_{\gamma,2})^{-1}
(\lambda_1B_{\gamma,1}+\lambda_2B_{\gamma,2}).
$$
It remains to note that the pairs $(C,\bA_{\gamma})$ and
$(C,\bA_{\gamma'})$ are not unitarily equivalent (which is shown by
another elementary calculation) unless $\gamma=\gamma'$.}
\end{example}

We conclude this section with characterizing Schur-class functions that
admit contractive commutative realizations.

              \begin{theorem}
              A Schur-class function  $S\in\cS_d(\cU,\cY)$ admits a contractive
              commutative realization if and only if it can be extended
              to a Schur-class function
              \begin{equation}
              \widehat{S}(\blam)=\begin{bmatrix}S(\blam) &
\widetilde{S}(\blam)\end{bmatrix}
              \in\cS_d(\cU\oplus\cF,\cY)
              \label{3.0n}
              \end{equation}
              such that the de Branges-Rovnyak space $\cH(K_{\widehat{S}})$ is
              ${\bf M}_{\blam}^*$-invariant and the inequality
              \begin{equation}
              \sum_{j=1}^d\|M_{\lam_j}^*f\|^2_{\cH(K_{\widehat{S}})}\le
              \|f\|^2_{\cH(K_{\widehat{S}})}-\|f(0)\|^2_{\cY}
              \label{3.0m}
              \end{equation}
              holds for every $f\in\cH(K_{\widehat{S}})$.
              \label{T:3.5a}
              \end{theorem}

              \begin{proof}
              Let $S$ admit a contractive commutative realization of the form
              \eqref{1.5a}.
              Extend the connecting operator $\bU$ of the form \eqref{1.7a} to a
              coisometric operator
              \begin{equation}
              \widehat{\bU}=\begin{bmatrix} A & B & \widetilde{B}\\ C & D
              & \widetilde{D}\end{bmatrix}\colon \; \begin{bmatrix} \cX \\
\cU\oplus\cF
              \end{bmatrix} \to \begin{bmatrix} \oplus_{1}^d\cX \\ \cY
\end{bmatrix}.
              \label{3.0s}
              \end{equation}
              The function
              \begin{equation}
              \widehat{S}(\blam)=\begin{bmatrix}D & \widetilde{D}\end{bmatrix}+
              C(I-Z(\blam)A)^{-1}Z(\blam)\begin{bmatrix}B &
\widetilde{B}\end{bmatrix}
              \label{3.0p}
              \end{equation}
              is an extension of $S$ in the sense of \eqref{3.0n}. The latter
realization is coisometric and commutative; thus ${\bf
M}_{\blam}^*$-invariance of $\cH(\widehat{S})$ and inequality
\eqref{3.0m} hold by Theorem \ref{T:comreal}.

Conversely, if $S$ can be extended to a Schur-class function $\widehat{S}$
with associated de Branges-Rovnyak space $\cH(K_{\widehat S})$
invariant under ${\bf M}_{\blam}^*$ and satisfying property
\eqref{3.0m}, we consider a weakly coisometric
commutative realization \eqref{3.0p} of $\widehat{S}$ (which exists by
Theorem \ref{T:comreal}) and restrict the input space to $\cU$. This gives
a contractive commutative realization for $S$.
\end{proof}

          \section{Realization for inner multipliers} \label{S:inner}

          In this section we focus on realization theory for inner
          multipliers.
            We first collect a couple of preliminary results needed for
the sequel.

	 \begin{theorem} \label{T:BBF1a}
	 Let $S \in {\mathcal S}_{d}(\cU, \cY)$, let $M_{S} \colon
             \cH_{\cU}(k_{d})
	 \to \cH_{\cY}(k_{d})$ be the multiplication operator defined in
	 \eqref{ms}, let $K_{S}$ denote the positive kernel given by
             \eqref{KS} and let
	   ${\mathbb M}_{S}$ denote the positive kernel
	   $$ {\mathbb M}_{S}(\blam, \bzeta) = \frac{1}{1 - \langle \blam,
	   \bzeta \rangle } \cdot S(\blam) S(\bzeta)^{*}.
	   $$
	   Then the reproducing kernel Hilbert spaces $\cH(K_{S})$ and
	   $\cH({\mathbb M}_{S})$ can be characterized as
	   $$ \cH(K_{S}) = \operatorname{Ran}\, (I - M_{S} M_{S}^{*})^{1/2},
	   \qquad \cH({\mathbb M}_{S}) = \operatorname{Ran}\, M_{S}
	   $$
	   with respective norms
	   \begin{align}
	     &  \| (I - M_{S} M_{S}^{*})^{1/2} f_{1} \|_{\cH(K_{S})} =
	     \| Q_{1} f_{1} \|_{\cH_{\cY}(k_{d})}\quad \text{for all} \quad
	     f_{1} \in\cH_{\cY}(k_{d}), \notag \\
	    & \| M_{S}f_{2}\|_{\cH({\mathbb M}_{S})} = \| Q_2
	    f_{2}\|_{\cH_{\cU}(k_{d})} \quad \text{for all} \quad
	    f_{2} \in\cH_{\cU}(k_{d})
	    \label{lifted-norm}
	    \end{align}
	    where $Q_{1}$ is the orthogonal projection onto
	    $(\operatorname{Ker}\, (I - M_{S}M_{S}^{*})^{1/2})^{\perp}$ and
	    $Q_{2}$ is the orthogonal projection onto $(\operatorname{Ker}\,
	    M_{S})^{\perp}$.

	   \end{theorem}

	 \begin{proof}
	The proof is based on a standard reproducing-kernel-space computation
	which we include for the sake of completeness.

	     Let $\cR_{D_{S^{*}}}$ denote the space
	    $ \operatorname{Ran}\, (I - M_{S} M_{S}^{*})^{1/2}$ with
	    norm given by the
	    $$
	     \| (I - M_{S} M_{S}^{*})^{1/2} f_{1}\|_{\cR_{D_{S^{*}}}} = \|
	     Q_{1} f_{1}\|.
	    $$
	    We also note the identity
	    $$ (I - M_{S} M_{S}^{*}) k_{d}(\cdot, \bzeta)y =
	    K_{S}(\cdot, \bzeta) y \text{ for all } \bzeta \in
	    {\mathbb B}^{d} \text{ and } y \in \cY.
	    $$
	    It follows that the set
	    $$
	   \cD := \overline{\operatorname{span}} \{(I - M_{S} M_{S}^{*})
	    k_{d}(\cdot, \bzeta) y \colon \bzeta \in {\mathbb
	    B}^{d},\, y \in \cY \}
	    $$
	    is dense in both $\cH(K_{S})$ and in $\cR_{D_{S^{*}}}$.
	    We then have, for all $\bzeta, \bzeta'
	    \in {\mathbb B}^{d}$ and $y, y' \in \cY$,
	    \begin{align*}
	   &	\langle (I - M_{S} M_{S}^{*}) k_{d}(\cdot, \bzeta) y,
		\,  (I - M_{S} M_{S}^{*}) k_{d}(\cdot, \bzeta') y'
		\rangle_{\cR_{D_{S^{*}}}} \\
		& \qquad =
		\langle (I - M_{S} M_{S}^{*}) k_{d}(\cdot, \bzeta) y, \,
		k_{d}(\cdot, \bzeta') y'  \rangle_{\cH_{\cY}(k_{d})}
		\\
		& \qquad = \langle K_{S}(\bzeta', \bzeta) y, y'
		\rangle_{\cY} \\
		& \qquad = \langle K_{S}(\cdot, \bzeta) y, K_{S}(\cdot,
		\bzeta') y' \rangle_{\cH(K_{S})} \\
		& \qquad = \langle (I-M_{S}M_{S}^{*}) k_{d}(\cdot, \bzeta)
		y,\, (I - M_{S} M_{S}^{*}) k_{d}(\cdot, \bzeta')y'
		\rangle_{\cH(K_{S})}.
	\end{align*}
	Hence the $\cR_{D_{S^{*}}}$ and the $\cH(K_{S})$ inner
	products agree on a common dense subset.  By taking closures
	and using completeness, it follows that $\cR_{D_{S^{*}}}$ and
	$\cH(K_{S})$ are equal to each other isometrically.

	The second statement in \eqref{lifted-norm} follows in a
	similar way by using the observation
	$$ M_{S} M_{S}^{*} k_{d}(\cdot, \bzeta) y = {\mathbb
	M}_{S}(\cdot, \bzeta) y.
	$$
	\end{proof}

A Schur-class function $S\in\cS_d(\cU,\cY)$ is said to be an {\em
inner multiplier} if the
multiplication operator $M_S$ (as an operator from $\cH_\cU(k_d)$ into
$\cH_\cY(k_d)$) is a partial isometry.

\begin{proposition}
\label{P:inner}
Let $S\in\cS_d(\cU,\cY)$. The following are equivalent:
          \begin{enumerate}
\item $S$ is inner.
\item $\cH(K_{S})$ is contained in $\cH_{\cY}(k_{d})$ isometrically.
\item $\cH(K_{S}) = ( \operatorname{Ran} \,M_{S})^{\perp}$ isometrically.
           \end{enumerate}
           In this case, $M_{S}M_{S}^{*}$ and  $I_{\cH_{\cY}(k_{d})} - M_{S}
M_{S}^{*}$ are the orthogonal projections onto
           the closed subspaces $\operatorname{Ran}\, M_{S}$ and
$(\operatorname{Ran} \,   M_{S})^{\perp}$ of $\cH_{\cY}(k_{d})$, respectively.
\end{proposition}

\begin{proof}  The multiplier $S$ being inner is equivalent to
        $\operatorname{Ran} M_{S}$ and
          $\operatorname{Ran}\, (I - M_{S} M_{S}^{*})^{1/2}$
          being closed
          subspaces of $\cH_{\cY}(k_{d})$ such that $M_{S}M_{S}^{*}$ and $I -
          M_{S}M_{S}^{*} = (I - M_{S} M_{S}^{*})^{1/2}$ are the orthogonal
          projections onto $\operatorname{Ran}\, M_{S}$ and
          $\operatorname{Ran}\, (I - M_{S} M_{S}^{*})$ respectively.  In
          this case the lifted-norm formulas \eqref{lifted-norm}
          lead to isometric inclusions
          of $\cH(K_{S})$ and of $\cH({\mathbb M}_{S})$ in
          $\cH_{\cY}(k_{d})$. The fact that $P_{\cH(K_{S})} = I -
          P_{\cH({\mathbb M}_{S})}$ tells us that $\cH(K_{S})$ and
          $\cH({\mathbb M}_{S})$ are orthogonal complements of each other.
         \end{proof}

          We are now ready for a realization characterization of inner
          multipliers.

          \begin{theorem}  \label{T:1.1a}
          An $\cL(\cU, \cY)$-valued function $S$ defined on ${\mathbb B}^{d}$
          is an inner multiplier if and only if it admits a weakly coisometric
          realization \eqref{1.5a} where:
          \begin{enumerate}
	  \item
       the $d$-tuple $\bA=(A_1,\ldots,A_d)$
          of the state space operators is commutative and is strongly stable
          (i.e., \eqref{stable} holds), and
          \item
          the output pair $(C, \bA)$ is isometric.
          \end{enumerate}
          \end{theorem}

          \begin{proof}
Suppose first that $S$ admits a realization \eqref{1.5a} with
$\bU = \sbm{A & B \\ C & D}$ weakly coisometric  with $\bA$ commutative
and strongly  stable  and  with \eqref{4.15} holding. By Proposition
\ref{R:1.3} we know that $K_{S}(\blam, \bzeta) = K_{C, \bA}(\blam,
\bzeta)$. Combining this equality with Theorem \ref{T:BBF1} (part (2)), we
conclude that the space $\cH(K_{S})= \cH(K_{C, \bA})$ is
included isometrically in $\cH_{\cY}(k_d)$. Therefore $S$ is
inner by Proposition \ref{P:inner}.

              Conversely, suppose that $S$ is inner.  Then, according to
              Proposition \ref{P:inner}, $\cH(K_{S})$ is
              isometrically equal to the orthogonal complement of
              $\operatorname{Ran}\, M_{S}$.  As $\operatorname{Ran}\, M_{S}$ is
              invariant under ${\mathbf M}_{\blam}$, it follows that
              $\cH(K_{S}) = ( \operatorname{Ran}\, M_{S})^{\perp}$ is ${\mathbf
              M}_{\blam}^{*}$-invariant.  Hence Theorem \ref{T:comreal}
              applies; we let $\bU = \sbm{A & B \\ C & D }$ be
              any weakly coisometric functional-model realization for $S$,
              that is with $\bA=(A_1,\ldots,A_d)$, $C$ and $D$ defined as in
\eqref{defCA}. Then $\bA$ is
              commutative since ${\mathbf M}_{\blam}^{*}$ is commutative.
              As has been already observed,
              $\cH(K_{S}) = ( \operatorname{Ran}\, M_{S})^{\perp}$ is
contained in
              $\cH_{\cY}(k_{d})$ isometrically.
              Therefore $\bA = {\mathbf M}_{\blam}^{*}|_{\cH(K_{S})}$ is
              strongly stable since ${\mathbf M}_{\blam}^{*}$ is strongly
              stable on $\cH_{\cY}(k_{d})$ by Proposition
              \ref{P:5.1} (part (3)). By part (2) in
              the same proposition, the pair $(G, {\mathbf M}_{\blam}^*)$ is
              isometric, i.e.,
\begin{equation}
I_{\cH_{\cY}(k_d)}-M_{\lambda_1}M_{\lambda_1}^*-\ldots
-M_{\lambda_d}M_{\lambda_d}^*=G^*G.
\label{gmisom}
\end{equation}
Since $G$ and $C$ are the operators of evaluation at the origin on
$\cH_{\cY}(k_d)$ and on $\cH(K_S)$ respectively, we have
$C=G\vert_{\cH(K_S)}$.
Then restricting operator equality \eqref{gmisom} to $\cH(K_S)$ we write
the obtained equality in terms of $\bA$ and $C$ as
$$
I_{\cH(K_S)}-A_1^*A_1-\ldots-A_d^*A_d=C^*C
$$
which means that the pair $(C,\bA)$ is isometric.
\end{proof}

The next theorem is a variant of Theorem \ref{T:3.5a} for the inner
case; the proof is much the same as that of Theorem \ref{T:3.5a} and
hence will be omitted.

\begin{theorem}  \label{T:.2}
       A Schur-class function $S \in \cS_{d}(\cU, \cY)$ admits a contractive,
       commutative realization of the form \eqref{1.5a} with $\bA =
       (A_{1}, \dots, A_{d})$ strongly stable and $(C, \bA)$ isometric if
       and only if $S$ can be extended to an inner multiplier
       $\widehat{S}(\blam)=\begin{bmatrix}S(\blam)
       & \widetilde{S}(\blam)\end{bmatrix}\in\cS_d(\cU\oplus\cF,\cY)$.
    \end{theorem}

If $S$ is inner, then, as we have seen in the proof of Theorem
\ref{T:1.1a}, any functional-model realization for $S$ yields a
commutative observable weakly coisometric realization for $S$.
We now show that any observable weakly coisometric realization for
$S$ necessarily is commutative.

\begin{theorem} \label{T:doesnot}
If $S \in \cS_{d}(\cU, \cY)$ is inner, then any observable weakly
coisometric realization of $S$ is  also commutative.
\end{theorem}

\begin{proof}
Let \eqref{1.5a} be an observable weakly coisometric realization of $S$.
Then $K_S=K_{C,\bA}$ (by Proposition \ref{R:1.3}) and therefore, since
$S$ is inner, the space $\cH(K_{C,\bA})$ is isometrically included into
$\cH_{\cY}(k_d)$. By Theorem \ref{T:3-1.2nc} (part (3)), the observability
operator ${\mathcal O}_{C,\bA}: \, \cX\to \cH_{\cY}(k_d)$ is a partial
isometry. Since the pair $(C,\bA)$ is observable, ${\mathcal O}_{C,\bA}$
is in fact an isometry.  Define the operators $T_1,\ldots,T_d$ on
$\cH(K_{C,\bA})$ and the operator $G: \, \cH(K_{C,\bA})\to \cY$ by
\begin{equation} \label{6.10a}
T_j {\mathcal O}_{C, \bA}x={\mathcal O}_{C,\bA}A_jx \; \; (j=1,\ldots,d),
\quad G{\mathcal O}_{C, \bA}x=Cx\quad       \text{for} \; \; x \in\cX.
          \end{equation}
Then for the generic element
$f={\mathcal O}_{C,\bA}x$ of $\cH(K_S)=\cH(K_{C,\bA})={\rm Ran} \,
{\mathcal O}_{C,\bA}$, we have
$$
f(\blam)=C(I-Z(\blam)A)^{-1}x,\qquad
f(0)=Cx=G{\mathcal O}_{C, \bA}x=Gf
$$
and therefore,
\begin{eqnarray*}
f(\blam)-f(0)&=&C(I-Z(\blam)A)^{-1}x-Cx\\
&=&C(I-Z(\blam)A)^{-1}Z(\blam)Ax \\
          &=&C(I-Z(\blam)A)^{-1}\sum_{j=1}^d \lambda_j A_jx \\
          &=&\sum_{j=1}^d  \lambda_j \cdot({\cO}_{C,\bA}A_jx)(\blam)\\
          &=&\sum_{j=1}^d  \lambda_j\cdot
          (T_j{\cO}_{C, \bA}x)(\blam)
          =\sum_{j=1}^d  \lambda_j\cdot (T_j f)(\blam)
          \end{eqnarray*}
which means that the $d$-tuple ${\bf T}=(T_1,\ldots,T_d)$ solves
the Gleason problem on $\cH(K_{C,\bA})$ and that $G$ is simply the
operator
of evaluation at the origin. Since the pair $(C,\bA)$ is contractive and
${\mathcal O}_{C,\bA}$ is isometric, it follows from \eqref{6.10a} that
the pair $(G,{\bf T})$ is also contractive. Now we recall a uniqueness
result from \cite{BBF1} (Theorem 3.22 there): {\em if ${\mathcal M}$ is a
backward-shift invariant subspace of $\cH_\cY(k_d)$ isometrically included
in $\cH_\cY(k_d)$, then the $d$-tuple ${\bf M}_{\blam}^*\vert_{\mathcal
M}=(M_{\lambda_1}^*\vert_{\mathcal
M},\ldots,M_{\lambda_d}^*\vert_{\mathcal
M})$ is the only contractive solution of the Gleason problem on
${\mathcal M}$}.
By this result applied to ${\mathcal M}=\cH(K_{C,\bA})=
\cH(K_{S})$ (which is backward-shift invariant and isometrically
included in $\cH_\cY(k_d)$ since $S$ is inner) we
conclude that $T_j=M_{\lambda_j}^*\vert_{\cH(K_{C,\bA})}$ for
$j=1,\ldots,d$. In particular, the tuple
${\bf T}$ is commutative and therefore the original state space tuple
${\bf A}$ is necessarily commutative.\end{proof}

We next observe that
in fact the {\em weakly coisometric} hypothesis in Theorem
\ref{T:doesnot} can be weakened to {\em contractive}.

\begin{corollary}  \label{C:doesnot}
If $S \in \cS_{d}(\cU, \cY)$ is inner, then any observable
contractive realization is commutative.
\end{corollary}

\begin{proof}
Assume that $S(\blam)=D + C (I - Z(\blam) A)^{-1} Z(\blam) B$ for
a contractive connecting operator  $\bU= \sbm{A & B \\ C & D}\colon \cX
\oplus \cU \to \cX^d \oplus \cY$. Extend $\bU$ to a
coisometric operator $ \widehat{\bU}$ as in \eqref{3.0s} and
consider its characteristic function $\widehat{S}$ (see
\eqref{3.0p}) which extends $S$ in the sense of \eqref{3.0n}. Since
$\widehat{\bU}$ is a contraction, $\widehat{S}(\blam)$ belongs to
$\cS_d(\cU\oplus\cF,\cY)$. By \eqref{3.0n},
$$
I_{\cY}-\widehat{S}(\blam)\widehat{S}(\blam)^*=I_{\cY}-S(\blam)S(\blam)^*-
\widetilde{S}(\blam)\widetilde{S}(\blam)^*\ge 0
$$
for almost all $\blam\in{\mathbb S}^d$. Since $S$ is inner, its boundary
values are coisometric almost everywhere on ${\mathbb S}^d$ (see
\cite{GRS}) and therefore, $\widetilde{S}(\blam)=0$ almost everywhere
on ${\mathbb S}^d$. Therefore, $\widetilde{S}\equiv 0$ and thus
$\widehat{S}$ is inner. The formula \eqref{3.0p} then gives an
observable coisometric (and therefore also weakly coisometric)
realization of the inner multiplier $\widehat S(\blam)$.
By Theorem \ref{T:doesnot} this realization necessarily is also
commutative, i.e., $\bA = (A_{1,} \dots, A_{d})$ is commutative.
Hence the original realization $\bU$ for $S(\blam)$ is commutative as
asserted.
\end{proof}
Theorems \ref{T:doesnot} and \ref{T:canonical-universal}
imply that, if $S\in\cS_d(\cU,\cY)$ is inner, then any contractive
observable realization of $S$ of the form \eqref{1.5a} is commutative
with operators $A_1,\ldots,A_d, C$  uniquely defined (up to
simultaneous unitary
equivalence) and with $D$ given by formulas \eqref{defCA}. The
nonuniqueness caused
by possible different choices of $B_1,\ldots,B_d: \, \cU\mapsto \cH(K_S)$
can be described explicitly. Let $S\in\cS_d(\cU,\cY)$ be inner, let $\bA =
(A_{1}, \dots, A_{d})$, $C$, $D$ be given as in \eqref{defCA}, let
${\mathcal D}$ be the subspace defined as in \eqref{canonical1}
(so that ${\mathcal D}^\perp:=\cH(K_S)^d\ominus{\mathcal D}$ is
characterized by \eqref{canonical2}), and let $B: \, \cU\to(\cH(K_S))^d$
be any operator so that $S$ can be realized in the form \eqref{1.5a}.
Then taking adjoints in \eqref{1.5a} gives
$$
B^*Z(\bzeta)^{*}(I - A^{*}Z(\bzeta)^{*})^{-1}C^{*}=S(\bzeta)^*-D^*
$$
which, on account of  \eqref{H-canonical}, can be written equivalently as
$$
B^*Z(\bzeta)^{*}K_S(\cdot, \, \bzeta)y=S(\bzeta)^*y-S(0)^*y\quad
(\bzeta\in{\mathbb B}^d, \; y\in\cY).
$$
Due to characterization \eqref{canonical1} of ${\mathcal D}$, the
latter formula completely determines the restriction of $B^*$ to
${\mathcal D}$:
\begin{equation}
B^*\vert_{\mathcal D}: \; Z(\bzeta)^{*}K_S(\cdot, \, \bzeta)y\to
S(\bzeta)^*y-S(0)^*y.
\label{c}
\end{equation}
Write $B^*: \, (\cH(K_S))^d\to \cU$ in the form
\begin{equation}
B^*=\begin{bmatrix}X & B^*\vert_{\mathcal D}\end{bmatrix}
\label{sep1}
\end{equation}
with $X=B^*\vert_{{\mathcal D}^\perp}: \; {\mathcal D}^\perp\to \cU$.
Next  we note the explicit formulas for the adjoints $A_j^{*}$'s
             \begin{equation}
             A_j^*={\mathcal P}_{\cH(K_S)}M_{\lam_j}\vert_{\cH(K_S)}\quad
                  (j=1,\ldots,d)
             \label{4.4}
             \end{equation}
(where ${\mathcal P}_{\cH(K_S)}$ stands for the orthogonal projection of
$\cH_{\cY}(k_d)$ onto $\cH(K_S)$) which are not available in the case of
general (noninner) Schur-class functions.
Indeed, since $\cH(K_S)$ is isometrically included in $\cH_{\cY}(k_d)$, we
have for every $h, \, g\in\cH(K_S)$,
             \begin{eqnarray*}
             \langle h, \; A_j^*g\rangle_{\cH(K_S)}&=&
             \langle A_jh, \; g\rangle_{\cH(K_S)}\\
             &=&\langle M_{\lam_j}^*h, \; g\rangle_{\cH(K_S)}\\
             &=&\langle M_{\lam_j}^*h, \; g\rangle_{\cH_{\cY}(k_d)}\\
             &=&\langle h, \; M_{\lam_j}g\rangle_{\cH_{\cY}(k_d)}\\
             &=&\langle h, \; {\mathcal P}_{\cH(K_S)}{M}_{\lam_j}g
                  \rangle_{\cH(K_S)}
             \end{eqnarray*}
                  and \eqref{4.4} follows. As a consequence of \eqref{4.4}
we get
             \begin{equation}
             A^*\vert_{\cD^{\perp}}=0.
             \label{4.5}
             \end{equation}
             Indeed, if $h=\sbm{ h_1 \\ \vdots \\ h_d}\in\cD^{\perp}$,
             it holds that $Z(\blam)h(\blam)\equiv 0$ (by the
                  characterization of $\cD^{\perp}$ in
                  \eqref{canonical2}) and then
             $$
             A^*h=\sum_{j=1}^d A_j^*h_j={\mathcal
P}_{\cH(K_S)}{M}_{\lam_j}h_j={\mathcal P}_{\cH(K_S)}(Zh)=0.
             $$
Now we define the operators $T_1: \, \cD \oplus \cY \to \cH(K_S)$ and
$T_{2} \colon \cD \oplus \cY \to \cU$ by
             \begin{equation}  \label{defT}
T_{1} = \begin{bmatrix} A^{*}|_{\cD} & C^{*}\end{bmatrix}\quad\mbox{and}
\quad T_{2} = \begin{bmatrix} B^*\vert_{\mathcal D} & S(0)^{*}
\end{bmatrix}
\end{equation}
and combining the two latter formulas with  \eqref{sep1} and \eqref{4.5},
we may write
the adjoint of the connecting operator $\bU=\sbm{A & B \\ C & D}$ as
\begin{equation}
\label{3.20}
\bU^*= \begin{bmatrix}  0 & T_{1} \\ X & T_{2}\end{bmatrix}
             \colon \begin{bmatrix} \cD^{\perp} \\ \cD
             \oplus \cY \end{bmatrix} \to \begin{bmatrix}\cH(K_S) \\ \cU
             \end{bmatrix}.
             \end{equation}
In the latter formula we have identified
$\begin{bmatrix} \cD^{\perp} \\ \cD
\oplus \cY \end{bmatrix}$ with $\begin{bmatrix}\cH(K_S)^{d}
\\
\cY \end{bmatrix}$. Every $X$ such that the matrix in \eqref{3.20} is
contractive leads to a contractive functional-model realization for $S$
(due to canonical choice \eqref{defCA} of $C$ and $\bA$)
which is automatically weakly coisometric. Therefore, the restriction of
$\bU^*$ to the space ${\mathcal D}\oplus\cY$, (that is, the operator
$\begin{bmatrix} T_1 \\ T_2\end{bmatrix}$) is isometric:
\begin{equation}
\label{sep2}
    T_1^*T_1+T_2^*T_2=I_{{\mathcal D}\oplus\cY}.
\end{equation}
Since the pair $(C,\bA)$ is isometric, it follows from \eqref{4.5}
and the formula for $T_1$ in \eqref{defT} that $T_1$ is coisometric:
\begin{equation}
\label{sep3}
    T_1T_1^*=A^*A+C^*C=I_{\cH(K_S)}.
\end{equation}
Then we also have
$T_1T_2^*T_2T_1^*=T_1(I-T_1^*T_1)T_1^*=I-I=0$, so that
\begin{equation}
\label{sep4}
T_1T_2^*=0.
\end{equation}
Now we invoke \eqref{3.20} and make use of \eqref{sep2}--\eqref{sep4}
to write the block-matrix formulas
\begin{equation}
I-\bU\bU^*=\begin{bmatrix} I-X^*X  & -X^*T_2 \\ -T_2^*X & 0
\end{bmatrix}
\label{sep5}
\end{equation}
and
\begin{equation}
I-\bU^*\bU=\begin{bmatrix}  0 & 0 \\ 0 & I-XX^*-T_{2}T_2^*\end{bmatrix}.
\label{sep6}
\end{equation}
 From the formula for $T_{2}$ in \eqref{defT} combined with the
          formula \eqref{c} for the action of $ B^*\vert_{\mathcal D}$ on a
generic generator of $\cD$, we see that
          $$
          \operatorname{\overline{Ran}}\, T_{2} =
\operatorname{\overline{span}}
          \{ S(\bzeta)^{*}y \colon \; \bzeta \in {\mathbb B}^{d},\, y \in
               \cY \}
          $$
          and hence
          \begin{equation} \label{KerT22}
\operatorname{Ker}\, T_{2}^{*} =(\operatorname{\overline{Ran}}
\, T_{2})^{\perp} = \{ u \in \cU \colon S(\blam) u \equiv 0 \}
=:\cU_{S}^{0}
          \end{equation}
Now it follows from \eqref{sep5} that $\bU^*$ of the
form \eqref{3.20} is contractive (isometric) if and only if $X$ is a
contraction (an isometry) from ${\mathcal  D}^\perp$ into (onto)
$\cU_S^0$. Then \eqref{sep6} implies that $\bU$ is unitary  if and only if
$X: \, {\mathcal  D}^\perp\to  \cU_S^0$ is unitary. The corresponding $B$
of the form \eqref{sep1} leads to weakly coisometric, coisometric and
unitary realizations for $S$.  We are led to the following result.
	   \begin{theorem} \label{T:4.3}
	   Let $S\in\cS_d(\cU,\cY)$ be inner, let $\bA = (A_{1},
	   \dots, A_{d})$, $C$, $D$ be given as in \eqref{defCA},
	   and let the subspaces
	   $\cU_S^{0}\subset\cU$ and $\cD$,
                  $\cD^\perp\subset\cH(K_S)^d$
	   be defined as in \eqref{KerT22}, \eqref{canonical1} and
	   \eqref{canonical2}. Then
	   \begin{enumerate}
             \item $S$ admits a coisometric functional-model
realization if and only if
$\dim \, \cD^{\perp}\le \dim \, \cU_S^0$.
\item $S$ admits a unitary functional-model realization if
and only if $\dim \, \cD^{\perp}= \dim \, \cU_S^0$.
	   \item $S$ admits a unique weakly coisometric
	   functional-model realization if and
	   only if $\cU_S^0=\{0\}$. In this case, the operator
	   $B\in\cL(\cU,\cH(K_S)^d)$ is defined by
$$
B^*\vert_{\mathcal D}: \;
Z(\bzeta)^{*}K_S(\cdot, \, \bzeta)y\to
S(\bzeta)^*y-S(0)^*y \quad \text{and}\quad
                  B^*|_{\cD^{\perp}}=0.
$$
This unique  weakly coisometric
             functional-model realization is never coisometric.
	   \end{enumerate}
	   \end{theorem}

\subsection{Beurling-Lax representation theorem for shift-invariant subspaces}

The Beurling-Lax theorem for the context of the Arveson space
$\cH_{\cY}(k_{d})$ asserts that any closed ${\mathbf M}_{\blam}$-invariant
subspace $\cM$ of $\cH_{\cY}(k_{d})$ can be represented in the form
\begin{equation}
{\mathcal M}=S\cdot \cH_{\cU}(k_d)
\label{sep8}
\end{equation}
for some inner multiplier
$S\in\cS_d(\cU,\cY)$ and an appropriately chosen coefficient space $\cU$
(see \cite{beurling, halmos, lax} for the classical case $d=1$ and
\cite{arv, mt, BBF1} for the case of general $d$).
We shall call any such $S$ a {\em  representer of ${\mathcal M}$}.  Here we
present a realization-theoretic proof of the
$\cH_{\cY}(k_{d})$-Beurling-Lax theorem as an application of Theorem
\ref{T:1.1a} (see \cite{BGR, BallRaney} for an illustration of this
approach for the case $d=1$).
We first need some preliminaries.

Suppose that $\bA = (A_{1}, \dots, A_{d})$ is a commutative
     $d$-tuple of bounded, linear operators on the Hilbert space $\cX$
     and that $(C, \bA)$ is an output stable pair.
     We define a {\em left-tangential functional
     calculus} $f \to (C^{*}f)^{\wedge
     L}(\bA^{*})$ on $\cH_{\cY}(k_{d})$ by
     \begin{equation}  \label{lefttanfunccal}
        (C^{*} f)^{\wedge L}(\bA^{*}) = \sum_{ \bn \in {\mathbb Z}^{d}_{+}}
        \bA^{* \bn} C^{*} f_{\bn}\quad \text{if } \quad f = \sum_{\bn \in
{\mathbb
        Z}^{d}_{+}} f_{\bn} \blam^{\bn} \in \cH_{\cY}(k_{d}).
   \end{equation}
   The computation
   \begin{align*}
       \left\langle \sum_{\bn \in {\mathbb Z}^{d}_{+}} \bA^{* \bn} C^{*}
   f_{\bn}, \; x \right \rangle_{\cX} & =
   \sum_{\bn \in {\mathbb Z}^{}_{+}} \left\langle f_{\bn}, \; C \bA^{\bn}
   x \right \rangle_{\cY} \\
   & = \sum_{\bn \in {\mathbb Z}^{d}_{+}} \frac{\bn !}{|\bn|!} \left\langle
    f_{\bn}, \frac{|\bn|!}{\bn !} C \bA^{\bn} x \right \rangle_{\cY} \\
    & = \langle f, \; \mathcal {O}_{C, \bA} x\rangle_{\cH_{\cY}(k_{d})}
   \end{align*}
   shows that the output-stability of the pair $(C, \bA)$ is exactly
   what is needed to verify that the infinite series in the definition
   \eqref{lefttanfunccal} of $(C^{*}f)^{\wedge L}(\bA^{*})$ converges
   in the weak topology on $\cX$. In fact the left-tangential
   evaluation with operator argument $f \to (C^{*}f)^{\wedge} (\bA^{*})$
   amounts to the adjoint of the observability operator:
   \begin{equation}
       (C^{*} f)^{\wedge L}(\bA^{*}) = (\cO_{C, \bA})^{*}
       f \quad\text{for}\quad f \in \cH_{\cY}(k_{d}).
   \end{equation}

   Given an output-stable pair $(C, \bA)$, define a subspace
   $\cM_{\bA^{*}, C^{*}} \subset  \cH_{\cY}(k_{d})$ by
   \begin{equation}  \label{MA*C*}
       \cM_{\bA^{*},C^{*}} = \{ f \in H_{\cY}(k_{d}) \colon
       (C^{*}f)^{\wedge L}(\bA^{*}) = 0 \}.
   \end{equation}
   An easy computation (using that $\bA$ is commutative) shows that
   $$
    (C^{*}[M_{\lam_{j}}f])^{\wedge L}(\bA^{*}) = A_{j}^{*} (C^{*}
    f)^{\wedge L}(\bA^{*}).
   $$
   Hence any subspace $\cM \subset \cH_{\cY}(k_{d})$ of the form
   $\cM = \cM_{\bA^{*}, C^{*}}$ as in \eqref{MA*C*} is
   ${\mathbf M}_{\blam}$-invariant.
   We now obtain the converse.

   \begin{theorem}  \label{T:shiftinv-homint} Suppose that $\cM$ is a
closed subspace of
   $\cH_{\cY}(k_{d})$ which is ${\mathbf M}_{\blam}$-invariant
   (i.e., $\cM$ is invariant under $M_{\lam_{j}} \colon f(\blam)
\mapsto \lam_{j} f(\blam)$
   for $j = 1, \dots, d$). Then there is a Hilbert space $\cX$, a
   commutative
   $d$-tuple of operators $\bA = (A_{1}, \dots, A_{d})$ on $\cX$ and
   an operator $C \colon \cX \to \cY$ so that
   \begin{enumerate}
       \item $\bA$ is commutative, i.e., $A_{i}A_{j} = A_{j} A_{i}$
       for $1 \le i,j \le d$,
       \item $\bA$ is strongly stable, i.e., $\bA$ satisfies \eqref{stable},
       and
       \item the subspace $\cM$ has the form $\cM_{\bA^{*}, C^{*}}$ as
       in \eqref{MA*C*}.
    \end{enumerate}
    Moreover, one choice of state space $\cX$ and operators $A_{j}
    \colon \cX \to \cX$ and $C \colon \cX \to \cY$ is
   \begin{equation}  \label{choice}
       \cX = \cM^{\perp}, \quad A_{j} =
       M_{\lam_{j}}^{*}|_{\cM^{\perp}} \quad\text{for} \; \; j = 1,
       \dots, d, \quad C \colon f \to f(0)\quad \text{for} \; \; f \in
       \cM^{\perp}.
   \end{equation}
    \end{theorem}

    \begin{proof}
    Define $\cX$, $\bA = (A_{1}, \dots, A_{d})$ and $C$ as in
    \eqref{choice}.  We note that $\cM^{\perp}$ is contained in
    $\cH_{\cY}(k_{d})$ isometrically and  that $C = G|_{\cM^{\perp}}$,
    $A_{j} = M_{\lam_{j}}^{*}|_{\cM^{\perp}}$ where $\cO_{G, {\mathbf
    M}_{\blam}^{*}}$ is the identity on $\cH_{\cY}(k_{d})$
   (see part (2) of Proposition \ref{P:5.1}).
    Hence in particular
    $\operatorname{Ran}\, \cO_{C, \bA} = \cM^{\perp}$.
    Taking orthogonal complements then gives
    $$
    \operatorname{Ker}\, (\cO_{C, \bA})^{*} = (\cM^{\perp})^{\perp} =
    \cM
    $$
    which in turn is equivalent to
    the characterization \eqref{MA*C*} for $\cM$.
   \end{proof}

   The following three examples illustrate how various  shift-invariant
   subspaces $\cM$ are characterized by some sort of concrete homogeneous
   interpolation conditions can be represented in the form
   $\cM_{\bA^{*},C^{*}}$; more general versions of these examples are
   discussed in \cite{BB-NYJ} in the context of
   nonhomogeneous interpolation problems of Nevanlinna-Pick type
   on domains in ${\mathbb C}^{d}$ of a more general form than
   ${\mathbb B}^{d}$.

   \begin{example} \label{E:hom1}
       {\em Let
   $${\bom}_{1} = (\omega_{1,1}, \dots,
       \omega_{1,d}), \dots, {\bom}_{n} = (\omega_{n,1}, \dots,
       \omega_{n,d}),
   $$
   be a collection of $n$ points in the unit ball ${\mathbb B}^{d}$,
   and let $x_{1}, \dots, x_{n}$ be a collection of linear functionals
   $x_{j} \colon \cY \to {\mathbb C}$ on $\cY$.
Define an associated shift-invariant subspace $\cM \subset
H^{2}_{\cY}$ by
$$ \cM = \{ f \in \cH_{\cY}(k_{d}) \colon x_{j}f({\bom}_{j})= 0
\text{ for } j = 1, \dots, n\}.
$$
Then it is easy to see that $\cM = \cM_{\bA^{*}, C^{*}}$ if one takes
$\cX = {\mathbb C}^{n}$ and
$$
   A_{j}^{*} = \begin{bmatrix} \omega_{1,j} & & \\ & \ddots & \\ & &
   \omega_{n,j} \end{bmatrix} \text{ for } j = 1, \dots, d, \qquad
   C^{*} = \begin{bmatrix} x_{1} \\ \vdots \\ x_{n} \end{bmatrix}.
$$}
\end{example}

\begin{example}  \label{E:hom2}
      Fix a point ${\bom} = (\omega_{1}, \dots,
      \omega_{d}) \in {\mathbb B}^{d}$ and let $x_{0}, \dots, x_{n-1}$
      be a collection of linear functionals $x_{j} \colon \cY \to
      {\mathbb C}$ $(j=0, \dots, n-1$) on $\cY$.  Associate a
      shift-invariant subspace $\cM \subset \cH_{\cY}(k_{d})$ by
      \begin{align*} \cM = &
	\{ f \in \cH_{\cY}(k_{d}) \colon
      x_{0}   \sum_{{\mathbf j} \colon |{\mathbf j}| = i}
      \frac{1}{i!} \frac{\partial^{i} f}{\partial \blam^{{\mathbf
      j}}} ({\bom})
      + x_{1}
      \sum_{{\mathbf j} \colon |{\mathbf j}| = i-1}
      \frac{1}{(i-1)!}
      \frac{\partial^{i-1}f}{\partial \blam^{\mathbf j}}
      ({\bom}) \\
      & \qquad + \cdots
      + x_{i} f({\bom})
       = 0 \text{ for
      } i=0,1, \dots, n-1\}.
      \end{align*}
      Then one can check that $\cM = \cM_{\bA^{*}, C^{*}}$ if one
      chooses $\cX = {\mathbb C}^{n}$ with
     $$
     A_{j}^{*} = \begin{bmatrix} \omega_{j} & & & \\ 1 & \omega_{j} & &
     \\ & \ddots & \ddots & \\ & & 1 & \omega_{j} \end{bmatrix} \text{
     for } j = 1, \dots, n, \qquad C^{*} = \begin{bmatrix} x_{0} \\
     x_{1} \\ \vdots \\ x_{n-1} \end{bmatrix}.
    $$
   \end{example}

   \begin{example} \label{E:hom3}
       {\em A more general example can be constructed as follows.
       Let ${\bom} = (\omega_{1}, \dots, \omega_{d})$ be a
       fixed point of ${\mathbb B}^{d}$, and
       let $E \subset {\mathbb Z}^{d}_{+}$ be a subset of indices which
       is {\em lower inclusive}, i.e.:  whenever $\bn \in E$ and $\bn -
       {\mathbf e}_{i} \in {\mathbb Z}^{d}_{+}$ (where ${\mathbf e}_{i} =
       (0, \dots, 1,
       \dots, 0)$ is the $i$-th
       standard basis vector for ${\mathbb R}^{d}$ with $i$-th
       component equal to $1$ and all other components equal to zero),
       then it is the case that also $\bn - {\mathbf e}_{i} \in E$.
       For each $\bn \in E$, let $x_{\bn} \colon \cY \to {\mathbb C}$
       be a linear functional on $\cY$. Define a polynomial $x(\blam)$
       with coefficients in $\cL(\cY, {\mathbb C})$ by
       $$
        x(\blam) = \sum_{\bn \colon \bn \in E} x_{\bn} (\blam -
        {\bom})^{\bn}.
       $$
       Define a subspace $\cM \subset
       H^{2}_{\cY}(k_{d})$ by
       $$
       \cM = \{ f \in \cH_{\cY}(k_{d}) \colon
       \frac{ \partial^{|\bn|}}{\partial \blam^{\bn}} \left\{ x(\blam) f(\blam)
       \right\}|_{\blam = {\bom}} = 0 \text{ for all } \bn \in E \}.
       $$
       Take the space $\cX$ to be equal to $\ell^{2}(E)$
       (complex-valued functions on the index set $E$ with
       norm-square-summable values) and define operators $A_{j}^{*} \in
       \cL(\ell^{2}(E))$ via square matrices with rows and columns
       indexed by $E$ as
       $$[ A_{j}^{*}]_{\bn, \bn'} = \delta_{\bn, \bn'} \omega_{j} +
       \delta_{\bn + {\mathbf e}_{j}, \bn'} \text{ for } j = 1, \dots,
       d
       $$
       and an operator $C^{*} \colon \cY \to \ell^{2}(E)$ as the column
       matrix
       $$  C^{*} = \text{col}_{\bn \in E} [x_{\bn}].
       $$
       Then it can be checked that $\cM = \cM_{\bA^{*}, C^{*}}$ for
       this choice of $(C, \bA)$.
       }
       \end{example}

   We now construct an inner multiplier solving a homogeneous
   interpolation problem via realization theory.

   \begin{theorem} \label{T:homint-BL}
       Suppose that $(C, \bA)$ is an isometric output-stable pair,
        with $\bA$ commutative and strongly stable.  Let $\cM =
       \cM_{\bA^{*},C^{*}} \subset \cH_{\cY}(k_{d})$ be given by
       \eqref{MA*C*}.  Then there is an input space $\cU$ and an inner
       Schur multiplier $S \in \cS_{d}(\cU, \cY)$
       so that $\cM = \operatorname{Ran} M_{S}$.  One such $S$ is given
       by
       $$
       S(\blam) = D + C (I- Z(\blam) A)^{-1} Z(\blam) B
       $$
where $A_1,\ldots,A_d$ and $C$ come from the given output pair $(C,
\bA)$ and
$B_1,\ldots,B_d$ are chosen so that the colligation
$\bU = \sbm{A & B \\ C & D }\colon \, \sbm{\cX \\ \cU}\to
\sbm{\cX^d \\ \cY}$ is weakly coisometric.
  In particular, one achieves a coisometric realization $\bU$
by choosing the input space $\cU$ and $\sbm{ B \\ D }$ so as to solve
   the Cholesky factorization problem:
   \begin{equation}  \label{Cholesky}
    \begin{bmatrix} B \\ D \end{bmatrix}
   \begin{bmatrix} B^{*} & D^{*} \end{bmatrix} = \begin{bmatrix}
       I_{\cX^{d}} & 0 \\ 0 & I_{\cY} \end{bmatrix} -
       \begin{bmatrix} A \\ C \end{bmatrix} \begin{bmatrix} A^{*} &
	 C^{*} \end{bmatrix}.
\end{equation}
\end{theorem}

    \begin{remark} {\em
        Note that the model output pair $(C, \bA)$ \eqref{choice}
        appearing in Theorem
        \ref{T:shiftinv-homint} is an isometric pair.
        In practice, however, one may given a subspace of the form
         $\cM_{\bA^{*},C}$ with $\bA$ commutative and strongly stable but
         without the pair
         $(C, \bA)$ being isometric.  If however it is the case that $(C,
         \bA)$ is {\em exactly observable} in the sense that the
         {\em observability gramian}
         $$
         {\mathcal G}_{C, \bA} = \cO_{C, \bA}^{*}\cO_{C, \bA} =
         \sum_{\bn \in {\mathbb Z}^{d}_{+}} \frac{
         |\bn|!}{\bn !} \bA^{\bn *} C^{*} C \bA^{\bn}
         $$
         is strictly positive
         definite, then the adjusted output pair $(\widetilde C, 
\widetilde \bA)$
         given by
         $$
         \widetilde A_{j} = H^{1/2} A_{j} H^{-1/2} \text{ for } j = 1,
         \dots, d, \qquad \widetilde C = C H^{-1/2} \text{ where } H: =
         {\mathcal G}_{C, \bA}
         $$
         is isometric and has all the other properties of the original output
         pair $(C, \bA)$, namely: $\widetilde \bA$ is strongly stable and
         $\cM = \cM_{\widetilde \bA^{*}, \widetilde C^{*}}$. Hence in
         practice the   requirement that the pair $(C, \bA)$ be
         isometric in
         Theorem \ref{T:homint-BL} can be replaced by the condition
         that $(C, \bA)$ is exactly observable. We note that in all the
         Examples \ref{E:hom1}, \ref{E:hom2} and \ref{E:hom3} above,
         the associated output pair $(C, \bA)$, while not isometric, is
         exactly observable.  A more complete
         discussion of this point can be found in \cite{BBF1}.}
\end{remark}

    \begin{proof}[Proof of Theorem \ref{T:homint-BL}]  Define $S(\blam)$ as in
    the statement of the theorem.  By Theorem \ref{T:1.1a}, $S$ is inner.
    By Proposition
    \ref{P:inner}, $(\operatorname{Ran} M_{S})^{\perp} = \cH(K_{S})$
    isometrically.  As $\bU$ is weakly coisometric, we also know that
    $\cH(K_{S}) = \cH(K_{C, \bA})$ by Proposition \ref{R:1.3}.  The
    space $\cH(K_{C, \bA})$ can in turn be identified as a set with
    $\operatorname{Ran}  \cO_{C, \bA}$ (see \cite[Theorem
    3.14]{BBF1}).  By hypothesis, $\cM = \operatorname{Ker}\, (\cO_{C,
    \bA})^{*}$; hence $\operatorname{Ran}\, \cO_{C, \bA} = \cM^{\perp}$.
    Putting all this together gives
    $(\operatorname{Ran} M_{S})^{\perp} = \cM^{\perp}$ and therefore
    $\operatorname{Ran} M_{S} = \cM$ as wanted.
    \end{proof}

    \begin{remark}  \label{R:inner}
{\em By the result of \cite{GRS}, it is known that inner
multipliers necessarily have nontangential boundary values on the unit
sphere ${\mathbb S}^{2d-1} = \partial {\mathbb B}^{d}$
which are almost everywhere (with respect to volume
measure on ${\mathbb S}^{2d-1}$) equal to partial isometries.
In the classical case $d=1$, inner multipliers are
characterized as those Schur-class functions whose
boundary values are partial isometries with a fixed
initial space (see \cite[Theorem C page 97]{RR}).  There appears to
be no analogous characterization in terms of boundary values of which
Schur-class functions are inner multipliers for the higher
dimensional case $d>1$; in short it is difficult to determine purely
from boundary-value behavior whether a given Schur-class multiplier
is an inner function or not.  On the other hand, Theorem
\ref{T:1.1a} enables to write down such inner functions
and Theorem \ref{T:homint-BL} enables us to write down the inner
representer for a given shift-invariant subspace.

It should be pointed out, however, that
the analogue of a Blaschke factor on the ball as the representer for
a codimension-1 ${\mathbf M}_{\blam}$-invariant subspace of
$\cH(k_{d})$ has been known for some
time  (see \cite{AlpayKaptanoglu}).
Indeed, for $a = \begin{bmatrix} a_{1} & \cdots & a_{d}\end{bmatrix}$
a point of ${\mathbb B}^{d}$ (viewed as a row matrix), with a little
bit of algebra one can see that the formula
for the $1 \times d$ Blaschke factor vanishing at $a$
$$
  b_{a}(z) = (1 - a a^{*})^{1/2}(1 - z a^{*})^{-1}(z-a) (I_{d} -
  a^{*} a)^{-1/2}
$$
(where the variable $z = \begin{bmatrix} z_{1} & \dots & z_{d}
\end{bmatrix} \in {\mathbb B}^{d}$ is also viewed as a row matrix)
appearing in \cite{AlpayKaptanoglu} can be written in realization form
\eqref{1.5a} with connecting operator $\bU$ given by
$$
\bU = \begin{bmatrix} A & B \\ C & D \end{bmatrix} =
\begin{bmatrix} a^{*} & (1-a^{*}a)^{1/2} \\ (1 - aa^{*})^{1/2} & -a
     \end{bmatrix} \colon \begin{bmatrix} {\mathbb C} \\ {\mathbb
     C}^{d} \end{bmatrix} \to \begin{bmatrix} {\mathbb C}^{d} \\
     {\mathbb C} \end{bmatrix}.
$$
We thus see that this Blaschke factor fits the prescription of
Theorem \ref{T:homint-BL} with
$(\bA^{*}, C^{*}) = (a, (1 - a^{*}a)^{1/2})$  and with $B,D$ chosen to
solve the Cholesky factorization problem \eqref{Cholesky} with the
resulting connecting operator $\bU$ unitary.  These Blaschke factors also
play an important role as the characterization of automorphisms of
the ball mapping the origin to a given point (see \cite[Theorem
2.2.2]{Rudin} and the references there for further history).}
\end{remark}

We next show how our analysis can be used to give a description of
all Beurling-Lax representers for a given shift-invariant subspace of
$\cH_{\cY}(k_{d})$.

\begin{theorem}  \label{T:repr}
Let ${\mathcal M}$ be a closed  ${\bf M}_{\blam}$-invariant subspace
of $\cH_{\cY}(k_d)$, let ${\mathcal N}={\mathcal
M}^\perp=\cH_{\cY}(k_d)\ominus {\mathcal M}$ and let
$$
A_{j} = M_{\lam_{j}}^{*}\vert_{\mathcal N} \;\; (j=1,\ldots,d),\quad
C \colon \; f \to f(0) \quad (f\in{\mathcal N}).
$$
Let $\cD$ be
the subspace of ${\mathcal N}^d$  given by \eqref{domV0} and let
$$
T:=\begin{bmatrix} A^*\vert_{\cD} & C^*\end{bmatrix}: \; \cD\oplus
\cY\to{\mathcal N}.
$$
Then:
\begin{enumerate}
\item Given a Hilbert space $\cU$, there exists an inner multiplier
$S\in\cS_d(\cU,\cY)$
satisfying \eqref{sep8} if and only if
\begin{equation}
\dim \, \cU\ge \dim \, {\operatorname{Ran}}\, (I-T^*T)^{\frac{1}{2}}.
\label{new3}
\end{equation}
\item If \eqref{new3} is satisfied, then all $S\in\cS_d(\cU,\cY)$
for which \eqref{sep8} holds are described by the formula
\begin{equation}
S(\blam)=\begin{bmatrix} C(I-Z(\blam)A)^{-1} &
I_\cY\end{bmatrix}(I-T^*T)^{\frac{1}{2}}G^*
\label{new3a}
\end{equation}
where $G$ is an isometry from
${\operatorname{Ran}}\, (I-T^*T)^{\frac{1}{2}}$ onto
$\operatorname{Ran}\, G\subset\cU$.
\item If $\dim \, \cU=\dim \,
\operatorname{Ran}\,(I-T^*T)^{\frac{1}{2}}$, then the function
$S\in\cS_d(\cU,\cY)$ such that \eqref{sep8} holds is defined uniquely up
to a constant unitary factor on the right.
\end{enumerate}
\end{theorem}

\begin{proof}  If $\cM$ is a closed ${\mathbf M}_{\blam}$-invariant
      subspace of $\cH_{\cY}(k_{d})$, then $\cM^{\perp}$ is
      isometrically included in $\cH_{\cY}(k_{d})$ and is isometrically
      equal to $\cH(K_{C,\bA})$ where $(C,\bA)$ are the model operators on
      $\cM^{\perp}$ as in \eqref{choice}.  On the other hand, as a
      consequence of Proposition \ref{P:inner}, we see that  the
      Schur multiplier $S$ is an inner-multiplier representer for $\cM$
      if and only if $\cH(K_{S}) = \cM^{\perp}$ isometrically.  Thus
      the problem of describing all inner-multiplier representers $S$
      for the given ${\mathbf M}_{\blam}$-invariant subspace $\cM$ is
      equivalent to the problem:  {\em describe all Schur-class multipliers
      $S$ such that $ \cH(K_{S}) = \cH(K_{C,\bA})$ isometrically}, where
$(C,\bA)$ is
      given by \eqref{choice}.  The various conclusions of Theorem
      \ref{T:repr} now follow as an application of Theorem 2.11
      from \cite{BBF2a} to the more special situation here (where $\bA$ is
      strongly stable and $\cH(K_{C,\bA})$ is contained in
      $\cH_{\cY}(k_{d})$ isometrically).
\end{proof}

\subsection{Examples}
In this section we illustrate a number of finer points concerning
realizations for inner functions with some examples.

\begin{example}  \label{E:4.4}
{\rm Here we give an example of an inner multiplier $S$ with a
unique, weakly coisometric, observable commutative realization $\bU$
which is not coisometric.  Thus {\em a fortiori} $S$ has no
observable realization which is also unitary. Let
	   \begin{equation}
S(\blam)=\begin{bmatrix} \lam_1^2 & \sqrt{2}\lam_1 \lam_2 &
              \lam_2^2\end{bmatrix}\quad\mbox{so that}\quad
K_S(\blam, \bzeta)=1+\lam_1\overline{\zeta}_1+
	   \lam_2\overline{\zeta}_2.
	   \label{4.1}
	   \end{equation}
	   Thus $S\in\cS_2(\C^3,\C)$ and the functions $\{1,\lam_1,
                  \lam_2\}$ form a basis for $\cH(K_S)$. It is readily seen
                  that $\cH(K_S)$ is invariant under the backward shifts
                  $M_{\lam_1}^*$ and $M_{\lam_2}^*$. Define matrices
	  \begin{align*}
	 &  A_1=\begin{bmatrix} 0 & 1 & 0 \\ 0 & 0 & 0 \\ 0 & 0 &
                  0\end{bmatrix},\quad
	   A_2=\begin{bmatrix} 0 & 0 & 1 \\ 0 & 0 & 0 \\ 0 & 0 &
                  0\end{bmatrix},\quad
	   B_1=\begin{bmatrix} 0 & 0 & 0 \\ 1 & 0 & 0 \\ 0 &
                  \frac{1}{\sqrt{2}} &
	   0\end{bmatrix}, \\
	   & B_2=\begin{bmatrix}0 & 0 & 0\\ 0 &
                  \frac{1}{\sqrt{2}}
	   & 0 \\ 0 & 0 & 1 \end{bmatrix}, \quad
	   C=\begin{bmatrix} 1 & 0 & 0\end{bmatrix}, \quad
	   D=\begin{bmatrix} 0 & 0 & 0\end{bmatrix}.
	  \end{align*}
	   Straightforward calculations show that
\begin{equation}
	   S(\blam)=D+C(I-\lam_1 A_1- \lam_2A_2)^{-1}(\lam_1 B_1+
                \lam_2B_2)
\label{july20}
\end{equation}
	   and that this realization is commutative, strongly stable
                  (in fact $A_1$ and $A_2$ are the matrices of operators
                  $M_{\lam_1}^*$ and $M_{\lam_2}^*$ restricted to 
$\cH(K_S)$ with
                  respect to the basis $\{1, \lam_1,\lam_2\}$ of
	   $\cH(K_S)$), observable and contractive (isometric, actually).
                  This realization is also weakly coisometric (by Proposition
                  \ref{R:1.3}), since
$$
C(I-\lam_1A_1- \lam_2A_2)^{-1}=\begin{bmatrix} 1 & \lam_1
& \lam_2\end{bmatrix}
$$
	   and then, \eqref{4.1} can be written in the form
$$
K_S(\blam, \bzeta)= C(I-\lam_1A_1- \lam_2A_2)^{-1}(I-\bar{\zeta}_1A^{*}_1-
\bar{\zeta}_2A^{*}_2)^{-1}C^*.
$$
In fact, this example amounts to the special case of Example
\ref{E:hom3} where ${\boldsymbol \omega} = (0,0)$,
$$
   E = \{ (0,0), (1,0), (0,1) \} \subset {\mathbb Z}^{2}_{+}
   $$
   and
   $$
     \cM = \left\{ f \in \cH(k_{2}) \colon f(0,0)=0, \,
     \frac{\partial f}{\partial \lam_{1}}(0,0) = 0, \, \frac{\partial
     f}{\partial \lam_{2}}(0,0) = 0 \right\}.
   $$
	   Thus $S$ is inner by Theorem \ref{T:1.1a}. Since $\cU_S^0=\{0\}$
	   (it is obvious), the above realization is the only
	   realization of $S$ (up to unitary
	   equivalence) which is weakly coisometric, commutative and
observable.
	   It is not coisometric and therefore $S$ does not admit  coisometric
	   (or {\em a fortiori} unitary) commutative observable realizations.

	   That $S$ cannot have a unitary realization can be seen directly as
	   follows.  If $S$ has a finite-dimensional unitary realization,
	   then necessarily $\dim \cX + \dim \cU = 2 \cdot \dim \cX + \dim
	   \cY$.  As $\dim \cU = 3$ and $\dim \cY = 1$, we then must have
	   $\dim \cX = 2$.  But this is impossible since $\dim \cH(K_{S}) =
	   3$.  On the other hand, a realization with $\dim \cX = \infty$
	   cannot be observable.}
	   \end{example}

	   \begin{example}   \label{E:4.5}
	   {\rm This example illustrates the nonuniqueness in the choice of
	   the input operator $B$ for a weakly coisometric, observable
	   realization $\bU$ for a given inner multiplier $S(\blam)$.

	   Let $S(\blam)$ be as in Example \ref{E:4.4} and let
	   $\widetilde{S}(\blam)=\begin{bmatrix} \lam_1^2 & \lam_1 \lam_2 &
	   \lam_2^2 & \lam_1 \lam_2 \end{bmatrix}$. Then
	   $$
	   K_{\widetilde{S}}(\blam, \bzeta)=1+\lam_1\overline{\zeta}_1+\lam_2
	   \overline{\zeta}_2=K_S(\blam, \bzeta)
	   $$
	   and thus the de Branges-Rovnyak spaces $\cH(K_{\widetilde{S}})$
                  and $\cH(K_S)$ coincide. We know from the previous example
                  that $\cH(K_S)$ is backward-shift-invariant and isometrically
                  included in $\cH(k_d)$; hence so also is
$\cH(K_{\widetilde{S}})$.
                  Therefore $\widetilde{S}$ is inner. Hence $\widetilde{S}$
                  admits  weakly coisometric commutative observable
	   realizations and, up to unitary similarity, the operators
                  $A_1$, $A_2$ and $C$ are as above while $D=\begin{bmatrix} 0 &
                  0 & 0 & 0\end{bmatrix}$.
	   Operators $B_1$ and $B_2$ are not uniquely defined. Note that
	   the space $\cU_{\widetilde{S}}^0$ is spanned by the vector
	   $u=\begin{bmatrix}0 & 1 & 0 & -1\end{bmatrix}^\top$. We already know
	   that $\cH(K_{\widetilde{S}})$ is $3$-dimensional and
                  $\{1,\lam_1, \lam_2\}$ is an orthonormal basis for
                  $\cH(K_{\widetilde{S}})$. We identify
	   $\cH(K_{\widetilde{S}})\oplus \cH(K_{\widetilde{S}})$ with
                  $\C^6$ upon identifying
	   the orthonormal basis
	   $$
	   \begin{bmatrix} 1 \\ 0\end{bmatrix}, \; \;
	   \begin{bmatrix} \lam_1 \\ 0\end{bmatrix}, \; \;
	   \begin{bmatrix} \lam_2 \\ 0\end{bmatrix}, \; \;
	   \begin{bmatrix} 0 \\ 1\end{bmatrix}, \; \;
	   \begin{bmatrix} 0 \\ \lam_1\end{bmatrix}, \; \;
	   \begin{bmatrix} 0 \\ \lam_2\end{bmatrix}
	   $$
	   for $\cH(K_{\widetilde{S}})\oplus \cH(K_{\widetilde{S}})$ with
                  the standard basis of $\C^6$. The subspace
	   $$
	   {\mathcal D}^\perp= \operatorname{span}
	   \left[\begin{array}{r} \lam_2 \\
                  -\lam_1\end{array}\right]\subset\cH(K_S)^2
	   $$
	   is identified with
	   the one dimensional subspace of $\C^6$ spanned by the vector
	   $x=\begin{bmatrix}0 & 0 & 1 & 0 & -1 & 0\end{bmatrix}^\top$. The
	   $5$-dimensional subspace ${\mathcal D}$ is then identified with
	   $$
	   \{x\}^\perp={\mathcal D}={\rm span} \,
	   \left\{\begin{bmatrix}\overline{\zeta}_1 &
	   \overline{\zeta}_1^2 &  \overline{\zeta}_1\overline{\zeta}_2 &
	   \overline{\zeta}_2 & \overline{\zeta}_1\overline{\zeta}_2 &
	   \overline{\zeta}_2^2\end{bmatrix}^\top: \; \zeta_1,
                  \zeta_2\in\C\right\}.
	   $$
	   The  operator $B^*=\begin{bmatrix} B_1^* & B_2^*\end{bmatrix}$ is
	   defined on $\{x\}^\perp$ by
	   $$
	   B^*\begin{bmatrix}\overline{\zeta}_1 &
	   \overline{\zeta}_1^2 &  \overline{\zeta}_1\overline{\zeta}_2 &
	   \overline{\zeta}_2 & \overline{\zeta}_1\overline{\zeta}_2 &
	   \overline{\zeta}_2^2\end{bmatrix}^\top=S(\zeta)^*-S(0)^*=
	   \begin{bmatrix}\overline{\zeta}_1^2 &
	   \overline{\zeta}_1\overline{\zeta}_2 &
	   \overline{\zeta}_2^2 &
                  \overline{\zeta}_1\overline{\zeta}_2\end{bmatrix}
	   $$
	   and $B^*$ must map $x$ into $\cU_{\widetilde{S}}$ contractively.
	   Since $\|x\|=\|u\| \, (=\sqrt{2})$, we set $B^*x=\alpha u$ with
	   $|\alpha|\le 1$ and this choice of $\alpha$ is the only
                  freedom we have. Thus, the matrices of $B_1$ and $B_2$ with
                  respect to the standard bases  are of the form
	   $$
	   B_1=\begin{bmatrix} 0 & 0 & 0 & 0 \\ 1 & 0 & 0 & 0 \\ 0 &
	   \frac{1+\alpha}{2} & 0 & \frac{1-\alpha}{2}
	   \end{bmatrix},\quad
	   B_2=\begin{bmatrix} 0 & 0 & 0 & 0 \\ 0 & \frac{1-\alpha}{2} & 0 &
	   \frac{1+\alpha}{2} \\ 0 & 0 & 1 & 0 \end{bmatrix}.
	   $$
	   If $|\alpha|=1$, then the corresponding realization is unitary.
	   Moreover, different $\alpha$'s lead to unitary realizations of
	   $\widetilde{S}$ that are not unitarily equivalent.}
	   \end{example}

\section{Characteristic functions of commutative row contractions}
\label{S:CharF}

In the operator model theory for commutative row contractions (see \cite{BES,
BES2}), one is given
a $d$-tuple of operators ${\mathbf T} = (T_{1}, \dots, T_{d})$ on a
Hilbert space $\cX$ for which the associated block-row matrix is
contractive:
$$ \|T\| \le 1 \text{ where } T  = \begin{bmatrix} T_{1} & \cdots &
T_{d} \end{bmatrix} \colon \cX^{d} \to \cX.
$$
Under certain conditions (that ${\mathbf T}$ be {\em completely
non-coisometric}---see
\cite{BES2, Popescu-var2}), the
associated characteristic function $\theta_{{\mathbf T}}(\blam)$ is a complete
unitary invariant for ${\mathbf T}$; recently extensions of the theory to still
more general settings have appeared (see \cite{BS, Popescu-var1,
Popescu-var2}) while the fully noncommutative setting is older (see
\cite{PopescuNF1, PopescuNF2}).  All this theory can be viewed as
multivariable analogues of the well-known, now classical operator
model theory of Sz.-Nagy-Foias \cite{NF}.  However, unlike the
fully developed theory in \cite{NF} for the classical case and unlike
the case for the fully noncommutative theory (see \cite{PopescuNF1,
PopescuNF2, Cuntz-scat}), none of the work for the multivariable
commutative setting provides a characterization of which Schur-class
functions arise as characteristic functions.

To define $\theta_{{\mathbf T}}$, we set
$\;  A = \begin{bmatrix} T_{1} & \ldots & T_{d} \end{bmatrix}^* \;$
and let ${\mathbf U} = \sbm{A & B \\ C & D } \colon \cX \oplus
{\mathcal D}_{T} \to \cX^{d} \oplus \cD_{T^{*}}$ be the {\em Halmos
unitary dilation} of $A$:
$$
B = D_{T}|_{{\mathcal D}_{T}} \colon \; {\mathcal D}_{T} \to
\cX^{d},\qquad  C  = D_{T^{*}} \colon \; \cX \to {\mathcal D}_{T^{*}},
$$
$$
D = -T|_{{\mathcal D}_{T}} \colon \; {\mathcal D}_{T} \to {\mathcal
       D}_{T^{*}},
$$
where
\begin{align*}
D_{T} &= (I_{\cX^{d}} - T^{*}T)^{1/2}, \qquad {\mathcal D}_{T} =
          \overline{\operatorname{Ran}}\, D_{T} \subset \cX^{d}, \\
D_{T^{*}} &= (I_{\cX} - T T^{*})^{1/2}, \qquad
       {\mathcal D}_{T^{*}} = \overline{\operatorname{Ran}}\, D_{T^{*}}
\subset \cX,
       \end{align*}
       and then $\theta_{{\mathbf T}}(\blam)$ is the transfer function 
associated
       with the colligation ${\mathbf U}$:
\begin{equation}
       \theta_{{\mathbf T}}(\blam) = [ -T + D_{T^{*}}(I - Z(\blam)T^{*})^{-1}
       Z(\blam) D_{T}]|_{{\mathcal D}_{T}} \colon {\mathcal D}_{T} \to
       {\mathcal D}_{T^{*}}.
\label{charf}
\end{equation}
       Since ${\mathbf U}$ is unitary, it follows that $\theta_{{\mathbf
T}}$ is in
       the Schur  class ${\mathcal S}_{d}({\mathcal
       D}_{T}, {\mathcal D}_{T^{*}})$.  More generally, a Schur-class
       function
       $S \in {\mathcal S}_{d}(\cU. \cY)$ is said to {\em coincide} with
       the characteristic function $\theta_{{\mathbf T}}(\blam)$ if there
are unitary
       identification operators
       $$ \alpha \colon {\mathcal D}_{T^{*}} \to \cY, \qquad \beta \colon
       {\mathcal D}_{T} \to \cU
       $$
       so that
       $$
       S(\blam) = \alpha \theta_{{\mathbf T}}(\blam) \beta^{*}.
       $$
       From our point of view, what is special about a Schur-class
       function $S
       \in {\mathcal S}_{d}(\cU, \cY)$ which coincides with a
       characteristic function $\theta_{{\mathbf T}}$ is that it is
required to have
       a {\em commutative unitary} realization.  An additional constraint
       follows from the fact that the unitary colligation in the
       construction of a characteristic function comes via the
       Halmos-dilation construction.  The following proposition summarizes
       the situation.  In general let us say that the Schur-class
       function
       $S \in {\mathcal S}_{d}(\cU, \cY)$ is {\em pure} if
       \begin{equation}  \label{pure}
           \|S(0)u\| = \| u \| \quad \text{for some} \quad u \in \cU\quad
           \Longrightarrow \quad u = 0.
       \end{equation}
       For the role of this notion in the characterization of
       characteristic functions for the classical case, see \cite[page
       188]{NF}.

       \begin{proposition} \label{P:charfunc}
	 A Schur-class function $S \in {\mathcal
	 S}_{d}(\cU, \cY)$ coincides with a characteristic function
	 $\theta_{{\mathbf T}}$ if and only if
\begin{enumerate}
\item $S$ has a realization
\eqref{1.5a} with ${\mathbf U}$ unitary and $\bA$ commutative, and
\item $S$ is pure, i.e., $S$ satisfies \eqref{pure}.
        \end{enumerate}
\end{proposition}

\begin{proof}
          We first note the following general fact:  {\em if ${\mathbf U} =
          \sbm{ A & B \\ C & D} \colon \cX \oplus \cU \to \cX^{d} \oplus
          \cY$ is unitary then the following are equivalent:
          \begin{enumerate}
	\item[(i)] $B$ is injective,
	\item[(ii)] $C^{*}$ is injective,
	\item[(iii)] $u \in \cU \;$ with $\; \|D u \| = \|u\| \;$ implies
             that $\; u = 0$.
           \end{enumerate}}
        To see this, note that the unitary property of ${\bf U}$ means that
         \begin{align}
         \begin{bmatrix} A^{*}A + C^{*}C & A^{*}B + C^{*}D \\
             B^{*}A + D^{*}C & B^{*}B + D^{*}D \end{bmatrix} & =
             \begin{bmatrix} I_{\cX} & 0 \\ 0 & I_{\cU} \end{bmatrix},
	   \notag \\
       \begin{bmatrix} A A^{*} + B B^{*} & A C^{*} + B D^{*} \\
            C A^{*} + D B^{*} & C C^{*} + D D^{*} \end{bmatrix} & =
            \begin{bmatrix} I_{\cX^{d}} & 0 \\ 0 & I_{\cY} \end{bmatrix}.
	  \label{relations}
        \end{align}
        From all these relations we read off
        \begin{align*}
         B u = 0 \quad &\Longrightarrow \quad \| D x \| = \| x\| \; \text{
and
         } \; C^{*}Dx= 0, \\
         \|D u \| = \| u \| \quad &\Longrightarrow \quad B u = 0 \; \text{
and
} \; C^{*}D u = 0, \\
         C^{*}y = 0 \quad &\Longrightarrow \quad  \|D^{*} y \| = \|y\| = \| D
(D^{*}y)\|  \;  \text{ and } \; B D^{*} y = 0.
        \end{align*}
        Hence any one of the conditions (i), (ii) or (iii) implies the
        remaining ones.

          Suppose now that $S \in {\mathcal S}_{d}(\cU, \cY)$ coincides
          with a characteristic function $\theta_{{\mathbf T}}$.  Then $S$ has a
          realization as in \eqref{1.5a} and \eqref{1.7a} for a
          connecting operator ${\mathbf
          U}$ of the form
          \begin{equation} \label{Halmos-col}
	  {\mathbf U} = \begin{bmatrix} A & B \\ C & D \end{bmatrix} =
	  \begin{bmatrix} I_{\cX^{d}} & 0 \\ 0 & \alpha \end{bmatrix}
	      \begin{bmatrix} T^{*} & D_{T} \\ D_{T^{*}} & -T \end{bmatrix}
		  \begin{bmatrix} I_{\cX} & 0 \\ 0 & \beta^{*}
		  \end{bmatrix}
           \end{equation}
           for a commutative row-contraction ${\mathbf T} = \{T_{1}, \dots,
           T_{d}\}$ where
            $\alpha \colon {\mathcal D}_{T^{*}} \to \cY$
             and $\beta \colon {\mathcal D}_{T} \to \cU$ are unitary, and
             where $\sbm{ T^{*} & D_{T} \\ D_{T^{*}} & -T}\colon \cX \oplus
             {\mathcal D}_{T} \to \cX^{d} \oplus {\mathcal D}_{T^{*}}$ is
             the Halmos dilation of $T^{*}$ discussed above.
             It is then obvious that ${\mathbf U}$ gives a commutative,
             unitary realization for $S$.  We also read off
             that any one (and hence all) of the conditions (i), (ii) and
             (iii) hold for ${\mathbf U}$.  As $D = S(0)$, the validity of
             condition (iii) implies that $S$ is pure.

             Conversely suppose that $S \in {\mathcal S}_{d}(\cU, \cY)$ has
             a commutative, unitary realization ${\mathbf U} = \sbm{ A & B
             \\ C & D } \colon \cX \oplus \cU \to \cX^{d} \oplus \cY$ and
             is pure.  As $D = S(0)$ and $S$ is pure, we read off that
             condition (iii) above holds, and hence also conditions (i) and
             (ii) hold for ${\mathbf U}$.  From the relations
             \eqref{relations} we have in particular
             $$ C^{*}C = I_{\cX} - A^{*} A, \qquad B B^{*} = I_{\cX^{d}} -
             A A^{*}.
             $$
             Hence we can define unitary operators
             $$ \alpha \colon {\mathcal D}_{A} \to
             \overline{\operatorname{Ran}}\, C = \cY, \qquad
             \beta \colon {\mathcal D}_{A^{*}} \to
             \overline{\operatorname{Ran}}\, B^{*} = \cU
             $$
             so that
             $$ \alpha D_{A} = C \quad\text{and}\quad \beta D_{A^{*}} = B^{*}.
             $$
             Then we also have
             $$ \alpha^{*} D \beta D_{A^{*}} = \alpha^{*} D B^{*} =
             -\alpha^{*} C A^{*} = -D_{A} A^{*} = -A^{*} D_{A^{*}}
             $$
             from which we get
             $$
               D = -\alpha A^{*} \beta^{*}.
            $$
            We conclude that ${\mathbf U}$ has the form \eqref{Halmos-col}
            with $T = (A_{1}^{*}, \dots, A_{d}^{*})$, and hence $S$
coincides with
            the characteristic function $\theta_{{\mathbf T}}$.
\end{proof}

         For the inner case, we can use the results on functional-model
         realizations obtained above to give a more intrinsic sufficient
         condition for a Schur-class function to be a characteristic function.

         \begin{theorem}  \label{T:charfunc-inner}
            Suppose that $S \in {\mathcal S}_{d}(\cU, \cY)$ is inner,
        $\dim \, {\mathcal D}^\perp=\dim \,
        \cU_S^0$ (where the subspaces $\cU_S^0\subset\cU$ and ${\mathcal
        D}^\perp\subset\cH(K_S)^d$  are defined in \eqref{KerT22} and
        \eqref{domV0}), and that $S$ is pure.  Then $S$ coincides with the
        characteristic function of a $*$-strongly stable, commutative
        row-contraction.
        \end{theorem}

        \begin{proof}  Given $S$ as in the hypotheses, we see from Theorem
            \ref{T:4.3} that $S$ has a functional-model realization
            ${\mathbf U} = \sbm{ A & B \\ C & D } \colon \cH(K_{S}) \oplus
            \cU \to \cH(K_{S})^{d} \oplus \cY$ such that ${\mathbf U}$ is
            unitary, $\bA$ is commutative and $\bA$ is strongly stable.
            Since by assumption $S$ is pure, we can apply Proposition
            \ref{P:charfunc} to conclude that $S$ coincides with
            $\theta_{{\mathbf T}}$, where ${\mathbf T} = (A_{1}^{*}, \dots,
A_{d}^{*})$.  As
            observed above, ${\mathbf A}$ is strongly stable, i.e.,
            ${\mathbf T}$ is $*$-strongly stable, and the theorem follows.
        \end{proof}


\begin{thebibliography}{99}

\bibitem{AM}
J.~Agler and J.E.~McCarthy,  \emph{Complete
Nevanlinna-Pick kernels}, J. Funct. Anal., {\bf  175} (2000),
no.1, 111--124.

\bibitem{AlpayKaptanoglu}  D.~Alpay and H.~Kaptano\u glu, {\em Some
finite-dimensional backward shift-invariant subspaces in the ball and
a related interpolation problem},  Integral Equations and Operator
Theory \textbf{42} no.~1 (2002), 1--21.


\bibitem{aron}
N.~Aronszajn, {\em Theory of reproducing kernels},
Trans.~Amer.~Math.~Soc., \textbf{68} (1950), 337--404.


\bibitem{arv}
W.~Arveson, {\em Subalgebras of $C^*$-algebras. III. Multivariable
operator theory}, Acta Math. \textbf{181} (1998), no. 2, 159--228.


\bibitem{BB-NYJ} J.A.~Ball and V.~Bolotnikov, {\em Nevanlinna-Pick
interpolation for Schur-AGler class functions on domains with matrix
polynomial defining function in ${\mathbb C}^{n}$}, New York J.~Math.
\textbf{11} (2005), 247-290.

\bibitem{BBF1}  J.A.~Ball, V.~Bolotnikov and Q.~Fang, {\em Multivariable
backward-shift invariant subspaces and observability operators},
Multidimensional Systems and Signal Processing, to appear.

\bibitem{BBF2a}  J.A.~Ball, V.~Bolotnikov and Q.~Fang, {\em
Transfer-function realization for multipliers of the Arveson space},
preprint.



\bibitem{BBF3} J.A.~Ball, V.~Bolotnikov and Q.~Fang, {\em Schur-class
multipliers on the Fock space: de Branges-Rovnyak reproducing kernel
spaces and transfer-function realizations}, in {\em Teberiu
Constantinescu Memorial Volume}, Theta, Bucharest, to appear.

\bibitem{BGR} J.A.~Ball, I.~Gohberg and L.~Rodman, {\em Interpolation
of Rational Matrix Functions}, \textbf{OT 45}, Birkh\"auser, Basel, 1990.


\bibitem{BallRaney} J.A.~Ball and M.~Raney, Discrete-time
dichotomous well-posed linear systems and
generalized Schur-Nevanlinna-Pick interpolation, {\em Complex
Analysis and Operator Theory}, to appear.

\bibitem{BTV}
J.A.~Ball, T.~T.~Trent and V.~Vinnikov, \emph{Interpolation and
commutant lifting for multipliers on reproducing kernels Hilbert
spaces}, in \emph {Operator Theory and Analysis} (Ed. H. Bart,
I. Gohberg and A.C.M. Ran), {\bf OT 122},
Birkh\"auser, Basel, 2001, pp. 89--138.


\bibitem{Cuntz-scat} J.A.~Ball and V.~Vinnikov, {\em Lax-Phillips
scattering and conservative linear systems: A Cuntz-algebra
multidimensional setting},  Mem. Amer. Math. Soc.  {\bf 178} (2005), no.
837.

\bibitem{beurling}
A. Beurling,
{\em On two problems concerning linear transformations in Hilbert space},
Acta Math. {\bf 81} (1949), 239--255.



\bibitem{BES} T.~Bhattacharyya, J.~Eschmeier and J.~Sarkar,
\emph{Characteristic function of a pure commuting contractive tuple},
Integral Equations and Operator Theory {\bf 53} (2005), no. 1,
23--32.

\bibitem{BES2} T.~Bhattacharyya, J.~Eschmeier and J.~Sarkar, \emph{On
c.n.c.~commuting contractive tuples}, arXiv:math.OA/0509162 v1, 7 Sep
2005.



\bibitem{BS} T.~Bhattacharyya and J.~Sarkar,
\emph{Characteristic function for polynomially contractive
commuting tuples}, J.~Math.~Anal.~Appl. \textbf{321} (2006), no. 1,
242-259.


\bibitem{dbr1}
{L. de} Branges and J.~Rovnyak,
Canonical models in quantum scattering theory,
in: {\em Perturbation Theory and its
Applications in Quantum Mechanics} (C.~Wilcox, ed.) pp. 295--392,
{Holt, Rinehart and Winston, New York}, 1966.

\bibitem{dbr2}
L.~de Branges and J.~Rovnyak,
{\em Square summable power series},
Holt, Rinehart and Winston, New York, 1966.



\bibitem{gleason}
A.~M.~Gleason, \emph{Finitely generated ideals in {B}anach algebras},
J. Math. Mech., \textbf{13} (1964), 125--132.

\bibitem{GRS}
D.C.~Greene, S.~Richter, and C.~Sundberg
\emph{The structure of inner multipliers on spaces with complete
Nevanlinna-Pick kernels}, J. Funct. Anal. {\bf 194} (2002), no. 2,
311--331.

\bibitem{halmos}
P.R. Halmos, {\em Shifts on Hilbert spaces},  J. Reine Angew. Math.
{\bf 208} (1961), 102--112.


\bibitem{lax}
P.D. Lax, {\em Translation invariant spaces},
Acta Math. {\bf 101} (1959), 163--178.



\bibitem{mt}
S.~McCullough and T.~T.~Trent, {\em Invariant subspaces and
Nevanlinna-Pick kernels}, J.~Funct.~Anal. {\bf 178} (2000), no. 1,
226--249.

\bibitem{NF} B.~Sz.-Nagy and C.~Foias, {\em Harmonic Analysis of
Operators on Hilbert Space}, North-Holland, Amsterdam-London, 1970.


\bibitem{PopescuNF1}
G.~Popescu, {\em Multi-analytic operators and some factorization
theorems}, Indiana Univ.~Math.~J. \textbf{38} (1989), no.~3,
6993-710.

\bibitem{PopescuNF2} G.~Popescu, {\em Characteristic functions for
infinite sequences of noncommuting operators}, J.~Operator
Theory \textbf{22} (1989) no. 1, 51-71.

\bibitem{Popescu-var1} G.~Popescu, \emph{Operator theory on noncommutative
varieties}, Indiana U.~Math.~J., to appear.

\bibitem{Popescu-var2} G.~Popescu, \emph{Operator theory on
noncommutative varieties II}, Proc. Amer. Math. Soc., to appear.

\bibitem{RR} M.~Rosenblum and J.~Rovnyak, {\em Hardy Classes and
Operator Theory}, Oxford University Press, Cambridge/New York, 1985.

\bibitem{Rudin} W.~Rudin, {\em Function Theory in the Unit Ball of
${\mathbb C}^{n}$}, Springer-Verlag, Berlin, 1980.



\end{thebibliography}
\end{document}